\documentclass[cmp]{svjour} 
\usepackage{amsmath,amssymb,
hyperref}
\usepackage{graphicx,color,subfig}

\def\makeheadbox{
}

\newenvironment{demo}[1]{\textit{Proof #1. }}{$\Box$}

\newtheorem{corol}[theorem]{Corollary}
\newtheorem{conj}[theorem]{Conjecture} 
\numberwithin{equation}{section}
\let\Diam\diamondsuit
\def\spin{\text{spin}}

\def\inj{\hookrightarrow}
\def\surj{\to\kern-.1em\llap{$\to$}}

\def\jrus{\leftarrow\kern-.1em\llap{$\leftarrow$}}
\def\strictsubset{\hbox{$\subseteq\kern-.4em\llap{${}_/$}$}}

\begin{document}
\title{{Discrete Riemann Surfaces and the Ising Model}}
\titlerunning{{Discrete Riemann Surfaces}}

\author{Christian \textsc{Mercat}}
\institute{Universit\'e Montpellier 2, France\\
\email{mercat@math.univ-montp2.fr}}

\date{Received 23 May 2000/ Accepted: 21 November 2000}
\renewcommand\communicated[1]{}
\communicated{M. E. {Fisher}}

\maketitle
\begin{abstract} We define a new theory of
discrete Riemann surfaces and present its basic results. The key idea
is to consider not only a cellular decomposition of a surface, but the
union with its dual. Discrete holomorphy is defined by a
straightforward discretisation of the Cauchy-Riemann equation. A lot
of classical results in Riemann theory have a discrete counterpart,
Hodge star, harmonicity, Hodge theorem, Weyl's lemma, Cauchy integral
formula, existence of holomorphic forms with prescribed holonomies.
Giving a geometrical meaning to the construction on a Riemann surface,
we define a notion of criticality on which we prove a continuous limit
theorem. We investigate its connection with criticality in the Ising
model. We set up a Dirac equation on a discrete universal spin
structure and we prove that the existence of a Dirac spinor is
equivalent to criticality.
\end{abstract}

\tableofcontents
\section{Introduction}
We present here a new theory of discrete analytic functions,
generalising to discrete Riemann surfaces the notion introduced by
Lelong-Ferrand~\cite{LF}.

Although the theory defined here may be applied wherever the usual
Riemann Surfaces theory can, it was primarily designed with
statistical mechanics, and particularly the Ising model, in
mind~\cite{McCW,ID}. Most of the results can be understood without any
prior knowledge in statistical mechanics. The other obvious fields of
application in two dimensions are electrical networks, elasticity
theory, thermodynamics and hydrodynamics, all fields in which
continuous Riemann surfaces theory gives wonderful results. The
relationship between the Ising model and holomorphy is almost as old
as the theory itself. The key connection to the Dirac equation goes
back to the work of Kaufman~\cite{K} and the results in this paper
should come as no surprise for workers in statistical mechanics; they
knew or suspected them for a long time, in one form or another. The 
aim
of this paper is therefore, from the statistical mechanics point of
view, to define a general theory as close to the continuous theory as
possible, in which claims as ``the Ising model near criticality
converges to a theory of Dirac spinors'' are given a precise meaning 
and a
proof, keeping in mind that such meanings and proofs already exist
elsewhere in other forms. The main new result in this context is that
there exists a discrete Dirac spinor near criticality in the finite
size Ising model,
\textbf{before} the thermodynamic limit is taken.
Self-duality, which enabled the first evaluations of the critical
temperature~\cite{KW,Ons,Wan50}, is equivalent to criticality at
finite size. It is given a meaning in terms of compatibility with
holomorphy.

The first idea in order to discretise surfaces is to consider
\textbf{cellular decompositions}.  Equipping a cellular decomposition
of a surface with a discrete metric, that is giving each edge a
length, is  sufficient if one only wants to do discrete harmonic
analysis. However it is not enough if one wants to define discrete
analytic geometry. The basic idea of this paper is to consider not
just the cellular decomposition but rather what we call its {\em
double,} {\it i.e.} the pair consisting of the cellular decomposition
together with its Poincar\'e dual. A discrete conformal structure is
then a class of metrics on the double where we retain only the ratio
of the lengths of dual edges\footnote{\rm By definition, a
discrete Riemann surface is a discrete surface equipped with a
discrete conformal structure in this sense.}. In Ising model terms, a
discrete conformal structure is nothing more than a set of interaction
constants on each edge separating neighbouring spins in an Ising 
model of
a given topology.

A function of the vertices of the double is said to be {\it discrete
holomorphic} if it satisfies the discrete Cauchy-Riemann equation, on
two dual edges $(x,x')$ and $(y,y')$, $$
\frac{f(y')-f(y)}{\ell(y,y')}=i\frac{f(x')-f(x)}{\ell(x,x')}.  $$
\begin{figure}[htbp]
\begin{center}\input{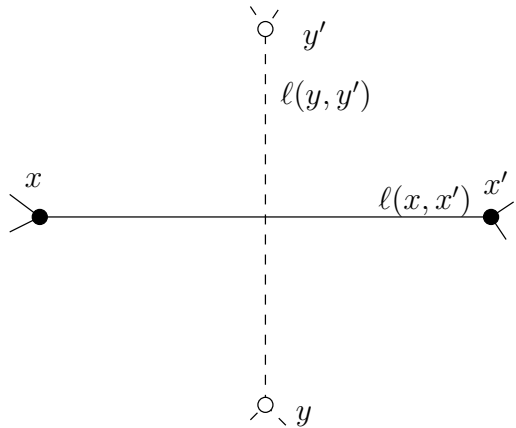}
\end{center}
\caption{The discrete Cauchy-Riemann equation.}         \label{fig:CR}

\end{figure}

This definition gives rise to a theory which is analogous to the
classical theory of Riemann surfaces. We define discrete differential
forms on the double, a Hodge star operator, discrete holomorphic
forms, and prove analogues of the Hodge decomposition and Weyl's
lemma. We extend to our situation the notion of pole of order one and
we prove existence theorems for meromorphic differentials with
prescribed poles and holonomies. Similarly, we define a Green
potential and a Cauchy integral formula.

Up to this point, the theory is purely combinatorial. In order to
relate the discrete and continuous theories on a Riemann surface, we
need to impose an extra condition on the discrete conformal structure
to give its parameters a geometrical meaning. We call this
semi-criticality in Sect.~\ref{sec:crit}. The main result here is
that the limit of a pointwise convergent sequence of discrete
holomorphic functions, on a refining sequence of semi-critical
cellular decompositions of the same Riemann surface, is a genuine
holomorphic function on the Riemann surface. If one imposes the
stronger condition of criticality on the discrete conformal structure,
one can define a wedge product between functions and $1$-forms which
is compatible with holomorphy.

Finally, for applications of this theory to statistical physics, one
needs to define a discrete analogue of spinor fields on Riemann
surfaces. In Sect.~\ref{sec:spin} we first define the notion of a
discrete spin structure on a discrete surface. It sheds an interesting
light onto the continuous notion, allowing us to redefine it in 
explicit
geometrical terms. In the case of a discrete Riemann surface we then
define a discrete Dirac equation, generalising an equation appearing
in the Ising model, and show that criticality of the discrete
conformal structure is equivalent to the existence of a local massless
Dirac spinor field.

In Sect.~\ref{sec:1st}, we present definitions and properties of the
theory which are purely combinatorial. First, in the empty boundary
case, we recall the definitions of dual cellular complexes, notions of
deRham cohomology. We define the double $\Lambda$, we present the
discrete Cauchy-Riemann equation, the discrete Hodge star on
$\Lambda$, the Laplacian and the Hodge decomposition. In
Subsect.~\ref{sec:DirNeu}, we prove Dirichlet and Neumann theorems,
the basic tools of discrete harmonic analysis.  In
Subsect.~\ref{sec:exist} we prove existence theorems for $1$-forms
with prescribed poles and holonomies.  In Subsect.~\ref{sec:wedge},
we deal with the basic difficulty of the theory: The Hodge star is
defined on $\Lambda$ while the wedge product is on another complex,
the diamond $\Diam$, obtained from $\Gamma$ or $\Gamma^*$ by the 
procedure of tile centering~\cite{GS87}. We prove Weyl's lemma, 
Green's identity and
Cauchy integral formulae.

In Sect.~\ref{sec:crit}, we define semi-criticality and criticality
and prove that it agrees with the usual notion for the Ising model on
the square and triangular lattices. 
We present Vorono\"\i~ and Delaunay semi-critical maps in order to
give examples and we prove the continuous limit theorem.  We prove 
that
every Riemann surface admits a critical map and give examples. On a
critical map, the product between functions and $1$-forms is
compatible with holomorphy and yields a polynomial ring, integration
and derivation of functions. We give an example showing where the
problems are.

In Sect.~\ref{sec:spin}, we set up the Dirac equation on discrete
spin structures. We motivate the discrete universal spin structure by
first showing the same construction in the continuous case. We show
discrete holomorphy property for Dirac spinors, we prove that
criticality is equivalent to the existence of local Dirac spinors and
present a continuous limit theorem for Dirac spinors.

\section{Discrete Harmonic and Holomorphic Functions}           
\label{sec:1st}
In this section, we are interested in properties of {\em
combinatorial} geometry. The constructions are considered up to
homeomorphisms, that is to say on a combinatorial surface, as opposed 
to
Sect.~\ref{sec:crit} where criticality implies that the discrete
geometry is embedded in a genuine Riemann surface.

\subsection{First definitions}
Let $\Sigma$ be an oriented surface without boundary.  A cellular
decomposition $\Gamma$ of $\Sigma$ is a partition of $\Sigma$ into
disjoint connected sets, called cells, of three types: a discrete set
of points, the vertices $\Gamma_0$; a set of non intersecting paths
between vertices, the edges $\Gamma_1$; and a set of topological discs
bounded by a finite number of edges and vertices, the oriented faces
$\Gamma_2$. A parametrisation of each cell is chosen, faces are mapped
to standard polygons of the euclidean plane, and edges to the segment
$(0,1)$; we recall particularly that for each edge, one of its two
possible orientations is chosen arbitrarily. We consider only locally
finite decompositions, {\it i.e.} any compact set intersects a finite
number of cells. In each dimension, we define the space of chains
$C_\bullet(\Gamma)$ as the $\mathbb{Z}$-module generated by the
cells. The boundary operator $\partial:C_k(\Gamma)\to C_{k-1}(\Gamma)$
partially encodes the incidence relations between cells. It fulfills
the boundary condition $\partial\partial=0$.

We now describe the \textbf{dual cellular decomposition} $\Gamma^*$ 
of a
cellular decomposition $\Gamma$ of a surface {\em without boundary}.
We refer to \cite{Veb} for the general definition.  Though we formally
use the parametrisation of each cell for the definition of the dual,
its combinatorics is intrinsically well defined.  To each face $F\in
\Gamma_2$ we define the vertex $F^*\in\Gamma^*_0$ inside the face $F$,
the preimage of the origin of the euclidean plane by the 
parametrisation
of the face.  Each edge $e\in\Gamma_1$, separates two faces, say $F_1,
F_2\in\Gamma_1$ (which may coincide), hence is identified with a
segment on the boundary of the standard polygon corresponding to
$F_1$, respectively $F_2$.  We define the dual edge $e^*\in\Gamma^*_1$
as the preimage of the two segments in these polygons, joining the
origin to the point of the boundary mapped to the middle of $e$. It is
a simple path lying in the faces $F_1$ and $F_2$, drawn between the
two vertices $F_1^*$ and $F_2^*$ (which may coincide), cutting no edge
but $e$, once and transversely. As the surface is oriented, to the
oriented edge $e$ we can associate the oriented dual edge $e^*$ such
that $(e,e^*)$ is direct at their crossing point. To each vertex
$v\in\Gamma_0$, with $v_1,\ldots,v_n\in\Gamma_0$ as neighbours, we
define the face $v^*\in\Gamma^*_2$ by its boundary $\partial
v^*=(v,v_1)^*+\ldots+(v,v_{k})^*+\ldots+(v,v_n)^*.$

\begin{remark} \rm

$\Gamma^*$ is a cellular decomposition of $\Sigma$~\cite{Veb}.
If we choose a parametrisation of the cells of $\Gamma^*$, we can
consider its dual $\Gamma^{**}$; it is a cellular decomposition
homeomorphic to $\Gamma$ but the orientation of the edges is
reversed. The bidual of $e\in\Gamma_1$ is the reversed edge
$e^{**}=-e$ (see Fig.~\ref{fig:duality}).
\begin{figure}[htbp]
\begin{center}\input{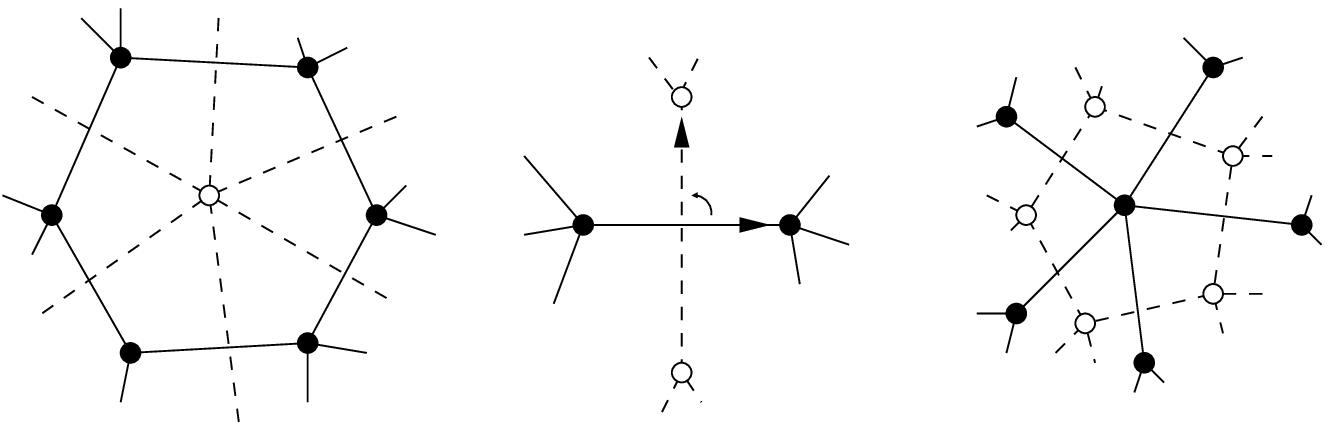}
\end{center}
\caption{Duality.}      \label{fig:duality}

\end{figure}
\end{remark}

The \textbf{double} $\Lambda$ of a cellular decomposition is the 
union of
these two dual cellular decompositions. We will speak of a $k$-cell of
$\Lambda$ as a $k$-cell of either $\Gamma$ or $\Gamma^*$.

A discrete \textbf{metric} $\ell$ is an assignment of a positive 
number
$\ell(e)$ to each edge $e\in\Lambda_1$, its \textbf{length}. For
convenience the edge with reversed orientation, $-e$, will be assigned
the same length: $\ell(-e):=\ell(e)$.  Two metrics
$\ell,\ell':\Lambda_1\to (0,+\infty)$ belong to the same discrete {\bf
  conformal structure} if the {\em ratio} of the lengths
$\rho(e):=\frac{\ell(e^*)}{\ell(e)}=\frac{\ell'(e^*)}{\ell'(e)}$, on
each pair of dual edges $e\in\Gamma_1, e^*\in\Gamma^*_1$ are equal.

A function $f$ on $\Lambda$ is a function defined on the vertices of
$\Gamma$ and of $\Gamma^*$. Such a function is said to be {\bf
  holomorphic} if, for every pair of dual edges $(x,x')\in\Gamma_1$ 
and
$(y,y')=(x,x')^*\in\Gamma^*_1$, it fulfills
$$ \frac{f(y')-f(y)}{\ell(y,y')}=i\frac{f(x')-f(x)}{\ell(x,x')}.
$$ It is the naive discretisation of the Cauchy-Riemann equation for a
function $f$, which is, in local orthonormal coordinates $(x,y)$:
$$\frac{\partial f}{\partial y}=i\frac{\partial f}{\partial x}.$$
Here, we understand two dual edges as being orthogonal.

This equation, though simple, was never considered in such a
generality. It was introduced by Lelong-Ferrand~\cite{LF} for the
decomposition of the plane by the standard square lattice
$\mathbb{Z}^2$. It is also called monodiffric functions; for
background on this topic, see~\cite{Duf}.  Polynomials of degree two,
restricted to the square lattice, give examples of monodiffric 
functions. See
also the works of Kenyon~\cite{Ken} and Schramm and
Benjamini~\cite{BS96} who considered more than lattices.

The usual notions of deRham cohomology are useful in this setup. We
said that $k$-chains are elements of the $\mathbb{Z}$-module
$C_k(\Lambda)$, generated by the $k$-cells, its dual space
$C^k(\Lambda):=\mathrm{Hom~}(C_k(\Lambda),\mathbb{C})$ is the space of
$k$-cochains. We will denote the coupling by the usual integral and
functional notations: $f(x)$ for a function $f\in C^0(\Lambda)$ on a
vertex $x\in\Lambda_0$; $\int_e\alpha$ for a $1$-form $\alpha\in
C^1(\Lambda)$ on an edge $e\in\Lambda_1$; and $\iint_F\omega$ for a
$2$-form $\omega\in C^2(\Lambda)$ on a face $F\in\Lambda_2$.

The \textbf{coboundary} $d:C^k(\Lambda)\to C^{k+1}(\Lambda)$ is 
defined
by the Stokes formula (with the same notations as before):
$$\int\limits_{(x,x')} df:= f\left(\partial(x,x')\right)=f(x')-f(x)
\qquad \iint\limits_F d\alpha:=\oint\limits_{\partial F}\alpha.
$$ As the boundary operator splits onto the two dual complexes
$\Gamma$ and $\Gamma^*$, the coboundary $d$ also respects the direct
sum $C^k(\Lambda)=C^k(\Gamma)\oplus C^k(\Gamma^*)$.

The Cauchy-Riemann equation can be written in the usual form
$*df=-idf$ for the following \textbf{Hodge star} $*:C^k(\Lambda)\to
C^{2-k}(\Lambda)$ defined by:
$$\int_e *\alpha:= -\rho(e^*)\int_{e^*}\alpha$$ We extend it to
functions and $2$-forms by:
$$\iint\limits_F *f:=f(F^*),\qquad *\omega
(x):=\iint\limits_{x^*}\omega.$$

As, by definition, for each edge $e\in\Lambda_1$,
$\rho(e)\rho(e^*)=\frac{\ell(e^*)}{\ell(e)}\frac{\ell(e)}{\ell(e^*)}=1$,

the Hodge star fulfills on $k$-forms,
$*^2=(-1)^{k(2-k)}\mathrm{Id}_{C^k(\Lambda)}$.

It decomposes $1$-forms into $-i$, respectively $+i$, eigenspaces,
called \textbf{type $(1,0)$}, resp. \textbf{type $(0,1)$ forms}:
$$C^1(\Lambda)=C^{(1,0)}(\Lambda)\oplus C^{(0,1)}(\Lambda).$$

The associated projections are denoted:
\begin{alignat}{2}
  \pi_{(1,0)}&=\frac12({\mathrm{Id}+i*})&:C^1(\Lambda)&\to
  C^{(1,0)}(\Lambda),\notag\\ 
  \pi_{(0,1)}&=\frac12({\mathrm{Id}-i*})&:C^1(\Lambda)&\to
  C^{(0,1)}(\Lambda).\notag
\end{alignat}

A $1$-form is \textbf{holomorphic} if it is closed and of type 
$(1,0)$:
$$\alpha\in\Omega^1(\Lambda)\iff d\alpha=0 \mathrm{~and~}
*\alpha=-i\alpha.$$ It is \textbf{meromorphic} with a pole at a vertex
$x\in\Lambda_0$ if it is of type $(1,0)$ and not closed on the face
$x^*$. Its \textbf{residue} at $x$ is defined by
$$\mathrm{Res}_x(\alpha):=\frac{1}{2i\pi}\oint_{\partial x^*}\alpha.$$
The residue theorem is merely a tautology in this context.

We define $d',d''$, the composition of the coboundary with the 
projections on
eigenspaces of $*$ as its holomorphic and anti-holomorphic parts:
\begin{alignat}{2}
  d'&:=\pi_{(1,0)}\circ d,&\qquad d''&:=\pi_{(0,1)}\circ d\notag\\ 
  \intertext{ from functions to $1$-forms,} d'&:=d\circ
  \pi_{(0,1)},&\qquad d''&:=d\circ \pi_{(1,0)}\notag
\end{alignat}
from $1$-forms to $2$-forms and $d'=d''=0$ on $2$-forms. They verify
${d'}^2=0$ and ${d''}^2=0$.

The usual discrete \textbf{laplacian}, which splits onto $\Gamma$ and
$\Gamma^*$ independently, reads $\Delta:=-d*d*-*d*d$ as expected. Its
formula for a function $f\in C^0(\Lambda)$ on a vertex
$x\in\Lambda_0$, with $x_1,\ldots,x_n$ as neighbours is
\begin{equation}
  \label{eq:lapldef}
  (\Delta f)(x)=\sum_{k=1}^n \rho(x,x_k) \left(f(x)-f(x_k)\right).
\end{equation}

As in the continuous case, it can be written in terms of $d'$ and 
$d''$
operators: For functions, $\Delta =i* (d'd''-d''d')$, in particular
holomorphic and anti-holomorphic functions are harmonic. The same 
result
holds for $1$-forms.

{\em In the compact case}, the operator $d^*=-*d*$ is the adjoint of
the coboundary with respect to the usual scalar product,
$(f,g):=\sum_{x\in\Lambda_0}f(x)\bar g(x)$ on functions, similarly on
$2$-forms and
$$ (\alpha,\beta):=\sum_{e\in\Lambda_1}\rho(e)\left(
  \int_e\alpha\right) \left( \int_e\bar\beta\right)
\mathrm{~on~}1\mathrm{-forms.}$$ It gives rise to the \textbf{Hodge
  decomposition},
\begin{proposition}[Hodge theorem]
  In the compact case, the $k$-forms are decomposed into orthogonal
  direct sums of exact, coexact and harmonic forms: $$C^k(\Lambda)=
  \mathrm{Im~}d\oplus^\perp\mathrm{Im~}d^*\oplus^\perp\mathrm{Ker~}\Delta,$$ 
harmonic forms are the closed and coclosed ones:
  $$\mathrm{Ker~}\Delta=\mathrm{Ker~}d\cap\mathrm{Ker~}d^*.$$ In
  particular the only harmonic functions are locally constant.
  Harmonic $1$-forms are also the sum of holomorphic and
  anti-holomorphic ones:
  $$\mathrm{Ker~}\Delta=\mathrm{Ker~}d'\oplus^\perp\mathrm{Ker~}d''.$$
\end{proposition}

Beware that $\Lambda$ being disconnected, the space of locally
constant functions is $2$-dimensional. The function $\varepsilon$ 
which
is $+1$ on $\Gamma$ and $-1$ on $\Gamma^*$ is chosen as the second
basis vector.

The proof is algebraic and the same as in the continuous case.  As the
Laplacian decomposes onto the two dual graphs, this result tells also
that for any harmonic $1$-form on $\Gamma$, there exists a unique
harmonic $1$-form on the dual graph $\Gamma^*$ such that the couple is
a holomorphic $1$-form on $\Lambda$, it's simply
$\alpha_{\Gamma^*}:=i*\alpha_\Gamma$.  These decompositions don't hold
in the non-compact case; there exist non-closed and/or non-co-closed,
harmonic $1$-forms.

\subsection{Dirichlet and Neumann problems} \label{sec:DirNeu}
\begin{proposition}[Dirichlet problem]
  Consider a finite connected graph $\Gamma$, equipped with a function
  $\rho$ on the edges, and  a certain non-empty set of points $D$
  marked as its boundary.  For any boundary function
  $f^\partial:(\partial\Gamma)_0\to\mathbb{C}$,
  there exists a unique function $f$, harmonic on $\Gamma_0\setminus
  D$ such that $f|_{\partial\Gamma}=f^\partial$. 
\end{proposition}

We refer to the usual laplacian defined by Eq.~\ref{eq:lapldef}.

If $f^\partial = 0$, the solution is the null function. Otherwise, it
is the minimum of the strictly convex, positive functional $f\mapsto
(df,df)$, proper on the non-empty affine subspace of functions which
agree with $f^\partial$ on the boundary.$\Box$

\begin{definition}            \label{defi:bord}
  Given $\Gamma$ a cellular decomposition of a compact surface with
  boundary $\Sigma$, define the double
  $\Sigma^2:=\Sigma\cup\bar\Sigma$, union with the opposite oriented
  surface, along their boundary. The double $\Gamma^2$ is a cellular
  decomposition of the compact surface $\Sigma^2$. Consider its dual
  $\Gamma^{2*}$ and define $\Gamma^*:=\Sigma\cap\Gamma^{2*}$. We don't
  take into account the faces of $\Gamma^{2*}$ which are not
  completely inside $\Sigma$ but we do consider the half-edges dual to
  boundary edges of $\Gamma$ as genuine edges noted
  $(\partial\Gamma^*)_1$ and define
  $(\partial\Gamma^*)_0:=\Gamma^{2*}_1\cap\partial\Sigma$ as the set
  of their boundary vertices.

A function $\rho$ on the edges of $\Gamma$ yields an extension to
$\Gamma^*_1$ by defining $\rho(e^*):=\frac 1{\rho(e)}$.
\end{definition}

\begin{remark} \rm
  $\Gamma^*$ is not a cellular decomposition of the surface; the
  half-edges dual to boundary edges do not bound any face of 
$\Gamma^*$.
\end{remark}

\begin{proposition}[Neumann problem]
  Consider $\Gamma$ a cellular decomposition of a disk, equipped with
  a function $\rho$ on its edges.  Choose a boundary vertex
  $y_0\in(\partial\Gamma^*)_0$, a value $f_{0}\in\mathbb{C}$, and on
  the set of boundary edges $e\in(\partial\Gamma^*)_1$, not incident
  to $y_0$, a $1$-form $\alpha$.

  Then there exists a unique function $f$, harmonic on
  $\Gamma^*\setminus(\partial\Gamma^*)_0$ such that $f(y_0)=f_{0}$ and
  $\int_e df=\int_e \alpha$ for all $e\in(\partial\Gamma^*)_1$ not
  incident to $y_0$.
\end{proposition}

It is a dual problem.  Let $e^*_0\in (\partial\Gamma^*)_1$, be the 
edge
incident to $y_0$ and $e_0\in (\partial\Gamma)_1$ its dual. Consider,
on the set of boundary edges $e\in (\partial\Gamma)_1$ different from
$e_0$, the $1$-form defined by $i*\alpha$. Integrating it along the
boundary, we get a function $f^\partial$ on $(\partial\Gamma)_0$, well
defined up to an additive constant. The Dirichlet theorem gives us a
function $f$ harmonic on $\Gamma_0\setminus(\partial\Gamma)_0$
corresponding to $f^\partial$.
Integrating the closed $1$-form $i*df$ on $\Gamma^*$  yields the
desired harmonic function $f$.$\Box$

\begin{remark} \rm
  The number of boundary points in $\Gamma$ is the same as in
  $\Gamma^*$, and as every harmonic function on $\Gamma$, when the
  surface is a disk, defines a harmonic function on $\Gamma^*$ such
  that their couple is holomorphic, unique up to an additive constant,
  the space of holomorphic functions, resp.  $1$-forms, on the double
  decomposition with boundary $\Lambda$ is $|(\partial 
\Lambda)_0|/2+1$,
  resp. $|(\partial \Lambda)_0|/2-1$ dimensional.

  The theorem is true for more general surfaces than a disk but the
  proof is different, see the author's PhD thesis~\cite{M}. There are
  $\ell^2$ versions of these theorems too.
\end{remark}



\subsection{Existence theorems}         \label{sec:exist}

We have very similar existence theorems to the ones in the continuous
case. We begin with the main difference:

\begin{proposition}                            \label{cor:2g}
  The space of discrete holomorphic $1$-forms on a compact surface 
without boundary
is of dimension twice the genus.
\end{proposition}

The Hodge theorem implies an isomorphism between the space of harmonic
forms and the cohomology group of $\Lambda$. It is the direct sum of
the cohomology groups of $\Gamma$ and of $\Gamma^*$ and each is
isomorphic to the cohomology group of the surface which is $2g$
dimensional on a genus $g$ surface. It splits in two isomorphic parts
under the type $(1,0)$ and type $(0,1)$ sum. As any holomorphic form
is harmonic, the dimension of the space of holomorphic $1$-forms is
$2g$.$\Box$

We can give explicit basis to this vector space as in the continuous
case~\cite{Sie}. To construct them, we begin with meromorphic forms:

\begin{proposition}                    \label{prop:exist1}
  Let $\Sigma$ be a compact surface with boundary.
 For each vertex $x\in\Lambda_0\setminus\partial\Diam$, and
  a simple path $\lambda$ on $\Lambda$ going from $x$ to the boundary
  there exists a pair of meromorphic $1$-forms $\alpha_x,\beta_x$ with
  a single pole at $x$, with residue $+1$ and which have pure
  imaginary, respectively real holonomies, along loops which don't 
have
  any edge dual to an edge of $\lambda$.
\end{proposition}

\begin{proposition}                    \label{prop:exist2}
  Let $\Sigma$ be a compact surface. For each pair of vertices
  $x,x'\in\Lambda_0$ with a simple path $\lambda$ on $\Lambda$ from
  $x$ to $x'$, there exists a unique pair of meromorphic $1$-forms
  $\alpha_{x,x'},\beta_{x,x'}$ with only poles at $x$ and $x'$, with
  residue $+1$ and $-1$ respectively, and which have pure imaginary,
  respectively real holonomies, along loops which don't have any edge
  dual to an edge of $\lambda$.
\end{proposition}

In both cases, the forms are $(\mathrm{Id}+i*)df$ with $f$ a solution
of a Dirichlet problem at $x$ (and $x'$) for $\alpha$ and of a Neumann
problem on the surface split open along the path $\lambda$ for
$\beta$. The uniqueness in the second proposition is given by the
difference: the poles cancel out and yield a holomorphic $1$-form with
pure imaginary, resp. real holonomies, so its real part, resp.
imaginary part, can be integrated into a harmonic, hence constant
function. So this part is in fact null. Being a holomorphic $1$-form,
the other part is null too.
We refer to the author's PhD thesis~\cite{M} for details.$\Box$

As in the continuous case, it allows us to construct holomorphic forms
with (no poles and) prescribed holonomies:

\begin{corol}                           \label{cor:PhiAB}
  Let $\mathcal{A},\mathcal{B}\in Z_1(\Lambda)$ be two 
non-intersecting simple 
  loops such that there exists exactly one edge of $\mathcal{A}$ dual
  to an edge of $\mathcal{B}$ (dual loops). There exists a unique
  holomorphic $1$-form $\Phi_{\mathcal{AB}}$ such that
  $\mathrm{Re}(\int_B\Phi_{\mathcal{AB}})=1$ and
  $\int_\gamma\Phi_{\mathcal{AB}}\in i\mathbb{R}$ for every loop
  $\gamma$ which doesn't have any edge dual to an edge of
  $\mathcal{A}$.
\end{corol}

We decompose $\mathcal{A}$ in two paths $\lambda_{x}^y$ and
$\lambda^{x}_y$. It gives us two $1$-forms $\beta_{x,y}$ and
$\beta_{y,x}$, then
\begin{equation}
  \Phi_{\mathcal{AB}}:=\frac1{2i\pi}(\beta_{x,y}+\beta_{y,x})
\end{equation}
fulfills the conditions.$\Box$
\subsection{The diamond  $\Diam$ and its  wedge product}      
\label{sec:wedge}
\begin{figure}[htbp]
\begin{center}\input{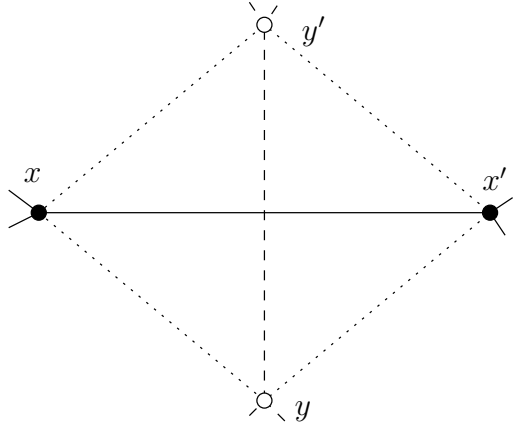}
\end{center}
\caption{The diamond $\Diam$.}  \label{fig:diamdefi}

\end{figure}

Following~\cite{Whit}, we define a wedge product, on another complex, 
the \textbf{diamond} $\Diam$,
constructed out of the double $\Lambda$:
Each pair of dual edges, say $(x,x')\in\Gamma_1$ and
$(y,y')=(x,x')^*\in\Gamma^*_1$, defines (up to homeomorphisms) a
four-sided polygon $(x,y,x',y')$ and all these constitute the faces of
a cellular complex called $\Diam$ (see Fig.~\ref{fig:diamdefi}).

On the other hand, from any cellular decomposition $\Diam$ of a
surface by four-sided polygons one can reconstruct the double
$\Lambda$. A difference is that $\Lambda$ may not be disconnected in
two dual pieces $\Gamma$ and $\Gamma^*$, it is so if each loop in
$\Diam$ is of even length; we will restrict ourselves to this simpler
case. This is not very restrictive because from a connected double,
refining each quadrilateral in four smaller quadrilaterals, one gets a
double disconnected in two dual pieces.

\begin{definition}
A \textbf{discrete surface with boundary} is defined by a 
quadrilateral
cellular decomposition $\Diam$ of an oriented surface with boundary
such that its double complex $\Lambda$ is disconnected in two dual 
parts. 
\end{definition}

This definition is a generalisation of the more natural previous
Definition~\ref{defi:bord}. It allows us to consider any subset of
faces of $\Diam$ as a domain yielding a discrete surface with
boundary. While any edge of $\Lambda$ has a dual edge, a vertex of
$\Lambda$ has a dual face if and only if it is an inner
vertex. Punctured surfaces can be understood in these terms too: An
inner vertex $v\in\Lambda_0$ is a puncture if it is declared as being
on the boundary and its dual face $v^*$ removed from $\Lambda_2$.

We construct a discrete wedge product, 
but while the Hodge star lives on the double $\Lambda$, the wedge
product is defined on the diamond $\Diam$: $\wedge:C^k(\Diam)\times
C^l(\Diam)\to C^{k+l}(\Diam).$ It is defined by the following
formulae, for $f,g\in C^0(\Diam)$, $\alpha,\beta\in C^1(\Diam)$ and
$\omega\in C^2(\Diam)$:
\begin{align*}
  (f\cdot g)(x):=&f(x)\cdot g(x)\qquad \mathrm{ ~for~} x\in\Diam_0,\\ 
  \int_{(x,y)}f\cdot\alpha:=& \frac{f(x)+f(y)}2
  \int{(x,y)}\alpha\qquad \mathrm{~for~} (x,y)\in\Diam_1,\\ 
  \iint\limits_{\hidewidth{(x_1,x_2,x_3,x_4)}\hidewidth}\alpha\wedge\beta
:=&\frac{1}{4}\sum_{k=1}^4
\int\limits_{{(x_{k-1},x_k)}}\alpha\;
\int\limits_{\hidewidth{(x_k,x_{k+1})}\hidewidth}\beta-
\int\limits_{{(x_{k+1},x_k)}}\alpha\;
\int\limits_{\hidewidth{(x_k,x_{k-1})}\hidewidth}\beta\\  
\iint\limits_{\hidewidth{(x_1,x_2,x_3,x_4)}\hidewidth}
  f\cdot\omega:=&\frac{\scriptstyle f(x_1)+f(x_2)+f(x_3)+f(x_4)}{4}
  \iint\limits_{\hidewidth{(x_1,x_2,x_3,x_4)}\hidewidth}\omega\\ 
  &\qquad \mathrm{~for~} (x_1,x_2,x_3,x_4)\in\Diam_2.
\end{align*}

\begin{lemma}
$d_\Diam$ is a derivation with respect to this
wedge product.
\end{lemma}

\begin{figure}[htbp]
\begin{center}\input{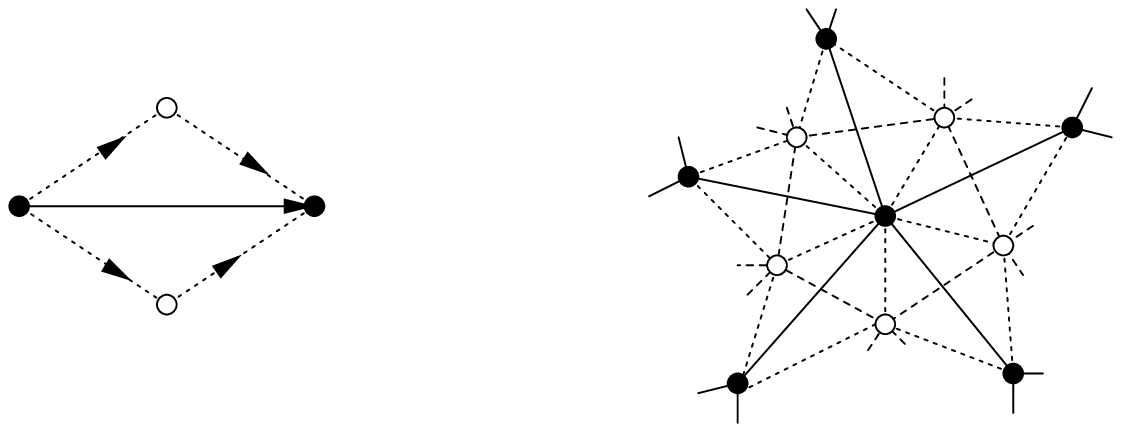}
\end{center}
\caption{Notations.}    \label{fig:avera}

\end{figure}

  To take advantage of this property, one has to relate
forms on $\Diam$ and forms on $\Lambda$ where the Hodge star is
defined. We construct an \textbf{averaging map} $A$ from 
$C^\bullet(\Diam)$
to $C^\bullet(\Lambda)$. The map is the identity for functions and
defined by the following formulae for $1$ and $2$-forms:
\begin{align}
  \int\limits_{{(x,x')}}A(\alpha_\Diam):= \frac12 \left(
    \int\limits_{{(x,y)}}+\int\limits_{\hidewidth{(y,x')}\hidewidth}
    +\int\limits_{\hidewidth{(x,y')}\hidewidth}+\int\limits_{{(y',x')}}
  \right) \alpha_\Diam, \label{def:avera1}\\ 
  \iint\limits_{\hidewidth{x^*}\hidewidth}A(\omega_\Diam):=
\frac12\sum_{k=1}^d\;\iint\limits_{{(x_k,y_k,x,y_{k-1})}}\omega_\Diam,                                                     
\label{def:avera2}
\end{align}
where notations are made clear in Fig.~\ref{fig:avera}. With this
definition, $d_\Lambda A=Ad_\Diam$, but the map $A$ is neither
injective nor always surjective, so we can neither define a Hodge star
on $\Diam$ nor a wedge product on $\Lambda$. An element of the
kernel of $A$ is given for example by $d_\Diam\varepsilon$, where
$\varepsilon$ is $+1$ on $\Gamma$ and $-1$ on $\Gamma^*$.

On the double $\Lambda$ itself, we have pointwise
multiplication between functions, functions and $2$-forms, and we 
construct an
{\em heterogeneous} wedge product for $1$-forms: with $\alpha,\beta\in
C^1(\Lambda)$, define $\alpha\wedge\beta\in C^1(\Diam)$ by
$$\iint\limits_{\hidewidth(x,y,x',y')\hidewidth}\alpha\wedge\beta:=
\int\limits_{\hidewidth(x,x')\hidewidth}\alpha
\int\limits_{\hidewidth(y,y')\hidewidth}\beta+
\int\limits_{\hidewidth(y,y')\hidewidth}\alpha
\int\limits_{\hidewidth(x',x)\hidewidth}\beta.$$

It verifies $A(\alpha_\Diam)\wedge A(\beta_\Diam)=\alpha_\Diam\wedge
\beta_\Diam$, the first wedge product being between $1$-forms on
$\Lambda$ and the second between forms on $\Diam$. Of course, we also
have for integrable $2$-forms:
$$\iint\limits_{\Diam_2}\omega_\Diam=\iint\limits_{\Gamma_2}A(\omega_\Diam)=
\iint\limits_{\Gamma^*_2}A(\omega_\Diam)=
\frac12\iint\limits_{\Lambda_2}A(\omega_\Diam).$$ And for a function
$f$,
$$\iint\limits_{\Diam_2}f\cdot\omega_\Diam=
\frac12\iint\limits_{\Lambda_2}A(f\cdot\omega_\Diam)=
\frac12\iint\limits_{\Lambda_2}f\cdot A(\omega_\Diam)$$ whenever
$f\cdot\omega_\Diam$ is integrable.

Explicit calculation shows that for a function $f\in C^0(\Lambda)$,
denoting by $\chi_x$ the characteristic function of a vertex
$x\in\Lambda_0$, $(\Delta f)(x)=-\iint_{\Lambda_2} f\cdot *\Delta
\chi_x.$ So by linearity one gets \textbf{Weyl's lemma}: a function 
$f$
is harmonic iff for any compactly supported function $g\in
C^0(\Lambda)$,
$$\iint\limits_{\Lambda_2} f\cdot \Delta g=0.$$

One checks also that the usual scalar product on compactly supported
forms on $\Lambda$ reads as expected:
$$ (\alpha,\beta):=\sum_{e\in\Lambda_1}\rho(e)\left(
  \int_e\alpha\right) \left(  
\int_e\bar\beta\right)=\iint\limits_{\Diam_2}\alpha\wedge*\bar\beta.$$

In some cases, for example, the decomposition of the plane by 
lattices,
the averaging map $A$ is surjective. We define the inverse map
$B:C^1(\Lambda)\to C^1(\Diam)/\mathrm{Ker~}A$ and $\Delta_\Diam
:=d_\Diam B*d$ and we then have
\begin{proposition}[Green's identity]
  For two functions $f,g$ on a compact domain $D\subset\Diam_2$,
  $$\iint_D (f\cdot \Delta_\Diam g -g\cdot \Delta_\Diam f) -
  \oint_{\partial D} (f\cdot B{*}d g-g\cdot B{*}d f)=0.
  $$
\end{proposition}

This means that for any representatives of the classes in
$C^1(\Diam)/\mathrm{Ker~}A$ the equality holds, but each integral
separately is not well defined on the classes.

\subsection{Cauchy integral formula}

\begin{proposition}
  Let $\Lambda$ a double map and $D$ a compact region of $\Diam_2$
  homeomorphic to a disc. Consider an interior edge $(x,y)\in D$;
  there exists a meromorphic $1$-form $\nu_{x,y}\in C^1(D\setminus
  (x,y))$ such that the holonomy $\oint_\gamma \nu_{x,y}$ along a
  cycle $\gamma$ in $D$ only depends on its homology class in
  $D\setminus (x,y)$, and $\oint_{\partial D} \nu_{x,y}=2i\pi$.
\end{proposition}

Consider the meromorphic $1$-form $\mu_{x,y}=\alpha_x+\alpha_y\in
C^1(\Lambda\cap D)$ defined by existence Theorem~\ref{prop:exist1} on
$D$. It is uniquely defined up to a global holomorphic form on $D$.
Its only poles are $x$ and $y$ of residue $+1$ so it verifies a 
similar
holonomy property, but on $\Lambda\cap D\setminus (x,y)$. We
define a $1$-form $\nu_{x,y}$ on $\Diam\cap D\setminus R$, such that
$\mu_{x,y}=A\nu_{x,y}$ in the following way: Let
$\int_{(x,{a})}\nu_{x,y}:=\lambda$, a fixed value, and for an edge
$(x',y')\in D_1$, with $x'\in\Gamma_0$, ${y'}\in\Gamma^*_0$, given two
paths in $D$, $\lambda_{x'}^x\in C_1(\Gamma)$ and $\lambda_{y}^{y'}\in
C_1(\Gamma^*)$ respectively from ${x'}$ to $x$ and from $y$ to ${y'}$,
$$\int_{({x'},{y'})}\nu_{x,y}:=
\int_{\lambda_{x'}^{x}}\mu_{x,y}+\int_{(x,A)}\nu_{x,y}+
\int_{\lambda_{y}^{y'}}\mu_{x,y}-\oint_{[\gamma]}\mu_{x,y},$$
where $[\gamma]$ is the homology class of
$\lambda_{x'}^{x}+(x,y)+\lambda_{y}^{y'}+({y'},{x'})$ on the punctured
domain.$\Box$

$\nu_{x,y}$ is the discrete analogue of $\frac {dz}{z-z_0}$ with
$z_0=(x,y)$. It is closed on every face of $D\setminus R$. By
definition, the average of $\nu_{x,y}$ on the double map is the
meromorphic form $A\nu_{x,y}=\mu_{x,y}$.

It allows us to state 
\begin{proposition}[Cauchy integral formula]
  Let $D$ be a compact connected subset of $\Diam_2$ and $(x,y)\in
  D_1$ two interior neighbours of $D$ with a non-empty boundary. For
  each function $f\in C^0(\Lambda)$,
\begin{equation}        \notag          \label{eq:cauchy}
  \oint_{\partial D} f\cdot\nu_{x,y}=\iint_D d''f\wedge \mu_{x,y}
  +2i\pi\frac{f(x)+f(y)}2.
\end{equation}
\end{proposition}

\begin{figure}[htbp]
\begin{center}\input{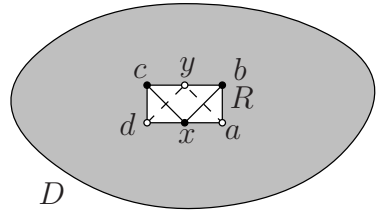}
\end{center}
\caption{The rectangle $R$ in a  domain $D$ defined by an edge
  $(x,y)\in\Diam_1$.}    
\label{fig:RABCD}
\end{figure}

The proof is straightforward: The edge $(x,y)$ bounds two faces in
$D$, let $R=(abcd)$ the rectangle made of these faces (see
Fig.~\ref{fig:RABCD}).

On $D\setminus R$,
$$d_\Diam (f\cdot\nu_{x,y})=d_\Diam f\wedge\nu_{x,y}+f\cdot
  d_\Diam\nu_{x,y}.$$

The $(1,0)$ part of $df$ disappears in the wedge product against the
holomorphic form $\mu_{x,y}$, so we can substitute
$$d_\Diam f\wedge\nu_{x,y} = d_\Lambda f\wedge A\nu_{x,y}=d''f\wedge 
\mu_{x,y}.$$
Integrating over $D$, as $\nu_{x,y}$ is closed on $D\setminus R$, we 
get:
$$\oint_{\partial D} f\cdot\nu_{x,y}=\iint_{D\setminus R} d''f\wedge
\mu_{x,y}+\oint_{\partial R} f\cdot\nu_{x,y}.$$

Explicit calculus shows that $\oint_{\partial R}
f\cdot\nu_{x,y}=\iint_Rd''f\wedge\mu_{x,y}+2i\pi\frac{f(x)+f(y)}2.$~$\Box$

\begin{remark} \rm
Since for all $\alpha\in C^1(\Diam)$, the locally constant function
$\varepsilon$ defined by $\varepsilon(\Gamma)=+1, \varepsilon(\Gamma^*)=-1$, verifies $\varepsilon\cdot\alpha=0$, an
integral formula will give the same result for a function $f$ and
$f+\lambda\varepsilon$. Therefore such a formula can not give access
to the value of the function at one point but only to its average
value at an edge of $\Diam$.
\end{remark}

\begin{corol}
  For $f\in\Omega(\Lambda)$ a holomorphic function, the Cauchy
  integral formula reads, with the same notations, 
$$\frac{f(x)+f(y)}2=\frac1{2i\pi}\oint_{\partial D} f\cdot\nu_{x,y}.$$
\end{corol}

The Green function on the lattices (rectangular, triangular,
hexagonal, Kagom\'e, square/octogon) is exactly known in terms of
hyperelliptic functions~(\cite{Hug} and references in Appendix~3). As
the potential is real, it means that the discrete Dirichlet problem on
these lattices can be exactly solved this way, once the boundary
values on the graph and its dual are given: if these values are real
and $\Gamma$ and imaginary on its dual, the solution is real on
$\Gamma$ and pure imaginary on the dual so the contributions $f(x)$
and $f(y)$ are simply the real and imaginary parts of the contour
summation respectively. Unfortunately, this pair of boundary values
are not independant but related by a Dirichlet to Neumann 
problem~\cite{CdV96}.
\section{Criticality}                                 \label{sec:crit}
The term criticality, as well as our motivation to investigate 
discrete
holomorphic functions, comes from statistical mechanics, namely the
Ising model. A critical temperature is defined that restrains the
interaction constants, interpreted here as lengths. We will see these
geometrical constraints in Sect.~\ref{sec:critcrit}.

Technically, as far as the continuous limit theorem is concerned, a
weaker property, called {\em semi-criticality} is sufficient, it gives
us a product between functions and forms. Moreover,
at criticality, this product will be compatible with holomorphy.

\subsection{Semi-criticality}
Define $C_\theta:=\{ (r,t):r\geq 0, t\in
\mathbb{R}/\theta\mathbb{Z}\}/(0,t)\sim(0,t')$ with the metric
$ds^2:=dr^2+r^2dt^2$ as the \textbf{standard cone of angle
  $\theta>0$}~\cite{Tro}.

The cones can be realized by cutting and pasting paper, demonstrating
their local isometry with the euclidean complex plane.

Let $\Sigma$ be a compact Riemann surface and $P\subset\Sigma$ a 
discrete
set of points. A \textbf{flat Riemannian metric with $P$ as conic
  singularities} is an atlas $\{ Z_{U_x}:U_x\to U'\subset
C_{\theta_x}>\}_{x\in P}$ of open sets $U_x\subset\Sigma$,
a neighbourhood of a singularity $x\in P$, into open sets of a 
standard
cone, such that the singularity is mapped to the vertex of the cone 
and
the changes of coordinates $C_{U,V}:U\cap V\to\mathbb{C}$ are 
euclidean
isometries.

There is a lot of freedom allowed in the choice of a flat metric for a
given closed Riemann surface $\Sigma$: Any finite set $P$ of points on
$\Sigma$ with a set of angles $\theta_x>0$ for every $x\in P$ such
that $2\pi\chi(\Sigma)=\sum_{x\in P} (2\pi - \theta_x)$, defines
uniquely a Riemannian flat metric on $\Sigma$ with these conic
singularities and angles~\cite{Tro}.

Consider such a flat riemannian metric on a compact Riemann surface
$\Sigma$ and $(\Lambda,\ell)$ a double cellular decomposition of
$\Sigma$ as before.
  \begin{definition}
    $(\Lambda,\ell)$ is a \textbf{semi-critical map} for this flat 
metric
    if the conic singularities are among the vertices of $\Lambda$ and
    $\Diam$ can be realized such that each face
    $(x,y,x',y')\in\Diam_2$ is mapped, by a local isometry $Z$
    preserving the orientation, to a four-sided polygon
    $(Z(x),Z(y),Z(x'),Z(y'))$ of the euclidean plane, the segments
    $[Z(x),Z(x')]$ and $[Z(y),Z(y')]$ being of lengths $\ell(x,x')$,
    $\ell(y,y')$ respectively and forming a direct orthogonal basis.
    We name $\delta(\Lambda,\ell)$ the supremum of the lengths of the
    edges of $\Diam$.
\end{definition}

  The local isometric maps $Z$ are discrete holomorphic.


 \textbf{Vorono\"\i~  and Delaunay complexes}~\cite{PS85} are
  interesting examples of semi-critical dual complexes. Any discrete
  set of points $Q$ on a flat Riemannian surface, containing the conic
  singularities, defines such a pair:

We first define two partitions $V$ and $D$ of $\Sigma$ into sets of
three types: $2$-sets, $1$-sets and $0$-sets, and then show that they
are in fact dual cellular complexes. They are defined by a real
positive function $m_Q$ on the surface, the {\em multiplicity}.

Consider a point $x\in\Sigma$; as the set $Q$ is discrete, the
distance $d(x,Q)$ is realized by geodesics of minimal length,
generically a single one. Let $m_Q(x)\in[1,\infty)$ be the number of 
such
geodesics.  If $m_Q(x)=1$, there exists a vertex $\pi(x)\in Q$ such
that the shortest geodesic from $x$ to $\pi(x)$ is the only geodesic
from $x$ to $Q$ with such a small length.

The Vorono\"\i~ $2$-set associated to a
vertex $v$ in $Q$, is $\pi^{-1}(v)$, that is to say the set of points 
of
$\Sigma$ closer to this vertex than to any other vertex in $Q$.
Each $2$-set of $V$ is a connected component of $m_Q^{-1}(1)$.

Likewise, the $1$-sets are the connected components of
$m_Q^{-1}(2)$. They are associated to pairs of points in $Q$.

The $0$-sets are the connected components of $m_Q^{-1}([3,+\infty))$.
Generically, they are associated to three points in $Q$.

$V$ is a cellular complex (see below) and the complex $D$ is its dual
(generically a triangulation), its vertices are the points in $Q$, its
edges are segments $(x,x')$ for $x,x'\in Q$ such that there exist
points equidistant and closer to them.

\begin{proposition} 
  The Vorono\"\i~ partition, on a closed Riemann surface with a flat
  metric with conic singularities, of a given discrete set of points
  $Q$ containing the conic singularities, is a cellular complex.
\end{proposition}

We have to prove that $2$-sets are homeomorphic to discs,
$1$-sets are segments and $0$-sets are points.

First,  $2$-sets are star-shaped, for every point $x$ closer to
$v\in Q$ than to any other point in $Q$, along a unique portion of a
geodesic, the whole segment $[x,v]$ has the same property.

$2$-sets are open, if $x$ is closer to $v\in Q$ than to any other
point in $Q$, as it is discrete, $d(x,Q\setminus v)-d(x,v)>0$. By
triangular inequality, every point in the open ball of this radius
centred at $x$ is closer to $v$ than to any other point in $Q$.

So $2$-sets are homeomorphic to discs.

Let $x$ be a point in a $1$-set. It is defined by exactly two 
portions of
geodesics $D,D'$ from $x$ to $y,y'\in Q$ (they may coincide). By
definition, the open sphere centred at $x$ containing $D\cup D'$
doesn't contain any point of $Q$ so it can be lifted to the universal
covering, where the usual rules of euclidean geometry tell us that the
set of points equidistant to $y$ and $y'$ around $x$ is a submanifold 
of
dimension $1$. As the surface is compact, if it is not a segment, it
can only be a circle. Then, it's easy to see that the surface is
homeomorphic to a $2$-sphere and that $y$ and $y'$ are the only points
in $Q$. But this is impossible because an euclidean metric on a
$2$-sphere has at least three conic singularities~\cite{Tro}.

 The same type of arguments shows that $0$-sets are isolated 
points.$\Box$

\begin{figure}[htbp]
\begin{center}\input{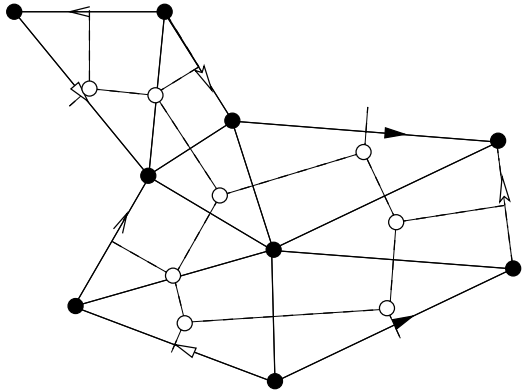}
\end{center}
\caption{The Vorono\"\i/Delaunay decompositions associated to two
  points on a genus two surface.} \label{fig:vorotore}

\end{figure}

\begin{proposition}
  Such Delaunay/Vorono\"\i~ dual complexes are semi-critical maps of 
the
  surface. Hence any Riemann surface admits semi-critical maps.
\end{proposition}

The edge in $V$ dual to $(x,x')\in D_1$ is a segment of their
mediatrix so is orthogonal to $(x,x')$. Hence, equipped
with the euclidean length on the edges, $(V,D)$ is a semi-critical 
map.$\Box$

\begin{remark} \rm
  Apart from Vorono\"\i/Delaunay maps, circle packings~\cite{CdV90}
  give another very large class of examples of interesting
  semi-critical decompositions (see Fig.~\ref{fig:circpack}).
\begin{figure}[htbp]
\begin{center}\input{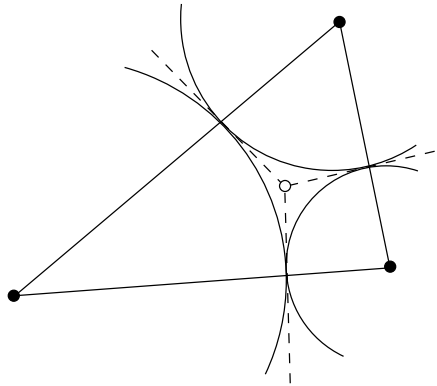}
\end{center}
\caption{Circle packing, the dual vertex to a face.}    
\label{fig:circpack}

\end{figure}

The semi-criticality of a double map gives a coherent system of angles
$\phi$ in $(0,\pi)$ on the oriented edges of $\Lambda$. An edge
$(x,x')\in\Lambda_1$ is the diagonal of a certain diamond;
$\phi(x,x')$ is the angle of that diamond at the vertex $x$. In
particular, $\phi(x,x')\not = \phi(x',x)$ {\it a priori}. They verify
that for every diamond, the sum of the angles on the four directions
of the two dual diagonals is $2\pi$ (see
Fig.~\ref{fig:systangsemi}). Then the conic angle at a vertex is
given by the sum of the angles over the incident edges.
\begin{figure}[htbp]
\begin{center}\input{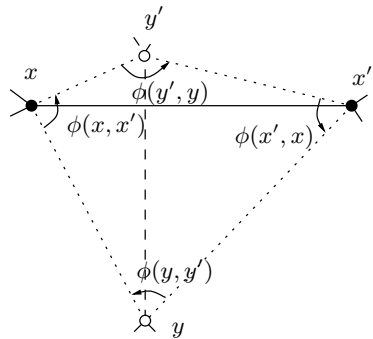}
\end{center}
\caption{A system of angles for a semi-critical map.}   
\label{fig:systangsemi}

\end{figure}
\end{remark}

\subsection{Continuous limit}
We state the main theorem, a converging sequence of discrete
holomorphic functions on a refining sequence of semi-critical maps of
the same Riemann surface, converges to a holomorphic function.
Precisely:

\begin{theorem}\label{th:lim}
  Let $\Sigma$ be a Riemann surface and
  $({}^k\!\Lambda,\ell_k)_{k\in\mathbb{N}}$ a sequence of
  semi-critical maps on it, with respect to the same flat metric with
  conic singularities. Assume that the lengths
  $\delta_k=\delta({}^k\!\Lambda)$ tend to zero and that the angles at
  the vertices of all the faces of the ${}^k\!\Diam$ are in the
  interval $[\eta, 2\pi-\eta]$ with $\eta>0$.
  
  Let $(f_k)_{k\in\mathbb{N}}$ be a sequence of discrete holomorphic
  functions $f_k\in\Omega({}^k\!\Lambda)$, such that there exists a
  function $f$ on $\Sigma$ which verifies, for every converging
  sequence $(x_k)_{k\in\mathbb{N}}$ of points of $\Sigma$ with each
  $x_k\in {}^k\!\Lambda_0$,
  $f\bigl(\lim_k(x_k)\bigr)=\lim_k\bigl(f_k(x_k)\bigr)$, then the
  function $f$ is \textbf{holomorphic} on $\Sigma$.
\end{theorem}

  Such a refining sequence is easy to produce (see
  Fig.~\ref{fig:refisemi}) but the theorem takes into account more
  general sequences. A more natural refining sequence, which mixes the
  two dual sequences is given by a series of tile centering
  procedures~\cite{GS87}: If one calls $\Diam/2$ the cellular
  decomposition constructed from $\Diam$ by replacing each tile by
  four smaller ones of half its size, and
  $\Gamma(\Diam/2),\Gamma^*(\Diam/2)$ the double cellular
  decomposition it defines, one has
  
$\Gamma(\Diam/2)=\Gamma(\Diam)\text{``}\cup\text{''}\Gamma^*(\Diam)$ 
and the
  interesting following sequence:
\begin{equation}
\begin{array}{ccccccccccccc}
\Gamma(\Diam)&\to&\Diam&\to&\Gamma(\Diam/2)&\to&\Diam/2&
\to&\cdots&\to&\Diam/2^n&\to&\dots\\
&\nearrow&&\searrow&&\nearrow&&
\searrow&&\nearrow&&\searrow&\\
\Gamma^*(\Diam)&&&&\Gamma^*(\Diam/2)&&&
&\cdots&&&&\cdots
\end{array}.
\end{equation}
The horizontal arrows correspond to tile centering procedures, and the
ascending, respectively descending arrows, to tile centering,
resp. edge centering procedures. This sequence is not that exciting
though since locally, the graph rapidly looks like a rectangular
lattice. More interesting inflation rules staying at criticality can
be considered too (see Fig.~\ref{fig:pentag2}).

The demonstration of the continuous limit theorem needs three lemmas:
\begin{figure}[htbp]
\begin{center}\input{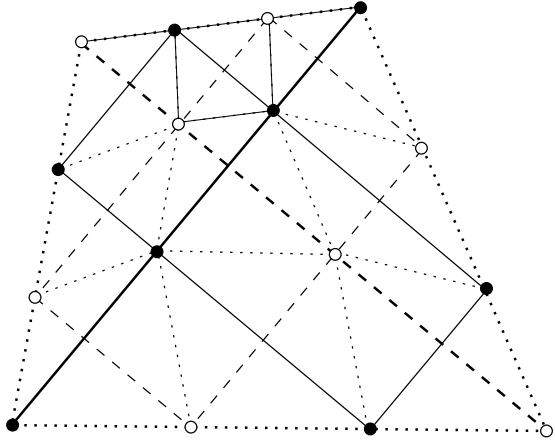}
\end{center}
\caption{Refining a semi-critical map.}
\label{fig:refisemi}

\end{figure}

\begin{lemma}   \label{lem:lim}
  Let $(f_k)_{k\in\mathbb{N}}$ be a sequence of functions on an open 
set
  $\Omega\subset\mathbb{C}$ such that there exists a function $f$ on
  $\Omega$ verifying, for every converging sequence
  $(x_k)_{k\in\mathbb{N}}$ of points of $\Omega$,
  $f\bigl(\lim_k(x_k)\bigr)=\lim_k\bigl(f_k(x_k)\bigr)$.  Then the
  function $f$ is continuous and uniform limit of $(f_k)$ on any
  compact.
\end{lemma}

Taking a constant sequence of points, we see that $(f_k)$ converges to
$f$ pointwise. So with the notations of the theorem, $(f_k(x_k))$
converges to $f(x)$ and $(f_l(x_k))_{l\in\mathbb{N}}$ to $f(x_k)$.
Combining the two, $(f(x_k))$ converges to $f(x)$ so $f$ is 
continuous.
If the convergence was not uniform on a compac sett, then there would
exist a converging sequence $(x_k)$ with $(f_k(x_k)-f(x_k))$ not
converging to zero. But $f$ is continuous in $x=\lim (x_k)$ and
$(f_k(x_k))$ converges to $f(x)$, which, combined, contradicts the
hypotheses.$\Box$

\begin{lemma}           \label{lem:diam}
  Let $(ABCD)$ be a four sided polygon of the Euclidean plane such 
that
  its diagonals are orthogonal and the vertices angles are in $[\eta,
  2\pi-\eta]$ with $\eta>0$. Let $(M,M')$ be a pair of points on the
  polygon. There exists a path on $(ABCD)$ from $M$ to $M'$ of minimal
  length $\ell$. Then
  $$\frac{MM'}\ell\geq\frac{\sin \eta}4.$$
\end{lemma}

\begin{figure}[htbp]
\begin{center}\input{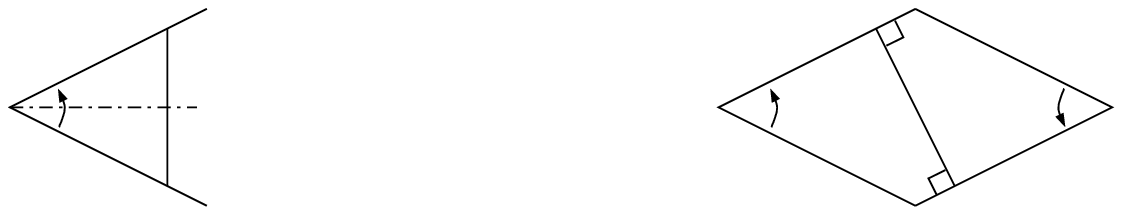}
\end{center}
\caption{The two extremal positions.}   \label{fig:diamlemma}

\end{figure}

It is a straightforward study of a several variables function. If the
two points are on the same side, $MM'=\ell$ and $\sin\eta\leq 1$. If
they are on adjacent sides, the extremal position with $MM'$ fixed is
when the triangle $MM'P$, with $P$ the vertex of $(ABCD)$ between
them, is isocel. The angle in $P$ being less than $\eta$,
$\frac{MM'}\ell\geq{\sin \frac\eta2}>\frac{\sin \eta}2.$ If the points
are on opposite sides, the extremal configuration is given by
Fig.~\ref{fig:diamlemma}.2., where $\frac{MM'}\ell=\frac{\sin
  \eta}4$.$\Box$

\begin{lemma}
  Let $(\Lambda,\ell)$ be any double cellular decomposition and
  $\alpha\in C^1(\Diam)$ a closed $1$-form. The $1$-form
  $f\cdot\alpha$ is closed for any holomorphic function
  $f\in\Omega(\Lambda)$ if and only if $\alpha$ is holomorphic.
\end{lemma}
Just check.

\begin{demo}{\ref{th:lim}}
  We interpolate each function $f_k$ from the discrete set of points
  ${}^k\!\Lambda_0$ to a function $\bar f_k$ of the whole surface,
  linearly on the edges of ${}^k\!\Diam$ and harmonicly in its faces.

  Let $(\zeta_k)$ be a converging sequence of points in $\Sigma$. Each
  $\zeta_k$ is in the adherence of a face of ${}^k\!\Diam$. Let
  $x_k,y_k$ be the minimum and maximum of $\mathrm{Re~} f_k$ around 
the
  face. By the maximum principle for the harmonic function
  $\mathrm{Re~}\bar f_k$,
  $$\mathrm{Re~} f_k(x_k)\leq \mathrm{Re~} \bar f_k(\zeta_k) \leq
  \mathrm{Re~} f_k(y_k).$$ Moreover, the distance between $x_k$ and
  $\zeta_k$ is at most $2\delta_k$, as well as for $y_k$. It implies
  that $(x_k)$ and $(y_k)$ converge to $x=\lim (\zeta_k)$,
  $(f_k(x_k))$ and $(f_k(y_k))$ to $f(x)$, and $(\mathrm{Re~} \bar
  f_k(\zeta_k))$ to $\mathrm{Re~} f(x)$; and similarly for its
  imaginary part. So, by Lemma~\ref{lem:lim}, the function $f$ is
  continuous, and is the uniform limit of $(\bar f_k)$ on every 
compact set. In
  particular, it is bounded on any compact.
  
  By the theorem of inessential singularities, since $f$ is continuous
  hence bounded on any compact set, and that conic singularities form 
a
  discrete set in $\Sigma$, to show that $f$ is holomorphic, we can
  restrict ourselves to each element $U\subset \Sigma$ of a euclidean
  atlas of the {\em punctured} surface (without conic
  singularities). We have an explicit coordinate $z$ on $U$.
   
  Let $\gamma$ be a homotopically trivial loop in $U$ of finite length
  $\ell$. We are going to prove that $\oint_\gamma fdz=0$. The theorem
  of Morera then states that $f$ is holomorphic.

  Let us fix the integer $k$.  By application of Lemma~\ref{lem:diam} 
on
  every face of ${}^k\!\Diam$ crossed by $\gamma$, we construct a loop
  $\gamma_k\in C_1({}^k\!\Diam)$, homotopic to $\gamma$, of length
  $\ell(\gamma_k)\leq \frac{4\ell}{\sin \eta}$ (see
  Fig.~\ref{fig:holgam}). As the diameter of a face of ${}^k\!\Diam$
  is at most $2\delta_k$, all these faces are contained in the tubular
  neighbourhood of $\gamma$ of diameter $4\delta_k$. Its area is
  $4\delta_k\ell$ and it contains the set $C$ of $\Sigma$ between
  $\gamma$ and $\gamma_k$.
\begin{figure}[htbp]
\begin{center}\input{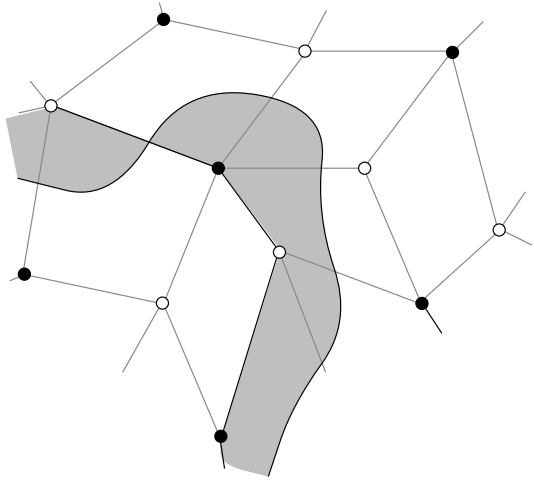}
\end{center}
\caption{The discretised path.}         \label{fig:holgam}

\end{figure}

Assume $f$ is of class $C^1$, on the compact $C$, $|\bar\partial f|$
is bounded by a number $M$. Applying Stockes formula to $fdz$,
$$|\oint_\gamma f(z)dz -\oint_{\gamma_k} f(z)dz|\leq
\iint_C|\bar\partial f(z)|dz\wedge d\bar z\leq M\times
4\delta_k\ell.$$ So $\oint_\gamma f(z)dz =\lim\oint_{\gamma_k}
f(z)dz$. Taking a sequence of class $C^1$ functions converging
uniformly to $f$ on $C$, we prove the same result for $f$ simply
continuous because all the paths into account are of bounded lengths.

As $(\bar f_k)$ converges uniformly to $f$ on $C$ and the paths are of
bounded lengths, we also have that $(|\oint_{\gamma_k} \bigl(\bar
f_k(z)- f(z)\bigr)dz|)_{k\in\mathbb{N}}$ tends to zero. But because
the interpolation is linear on edges of ${}^k\!\Diam$,
$\oint_{\gamma_k} \bar f_k(z)dz=\oint_{\gamma_k} f_kdZ$, the second
integral being the coupling between a $1$-chain and a $1$-cochain of
${}^k\!\Diam$.  But since $f_k$ and $dZ$ are discrete holomorphic,
$f_kdZ$ is a closed $1$-form, and $\oint_{\gamma_k} f_kdZ=0$. So
$\oint_{\gamma_k} f_k(z)dz$ tends to zero and
$$\oint_\gamma f(z)dz=0.$$
\end{demo}

\subsection{Criticality}                        \label{sec:critcrit}


\begin{proposition}            \label{th:holoz}
  Let $\alpha$ be a holomorphic $1$-form, $f\cdot\alpha$ is 
holomorphic
  for any holomorphic function if and only if
  $\int_{(y,x)}\alpha=\int_{(x',y')}\alpha$ for each pair of dual
  edges $(x,x'),(y,y')$.
\end{proposition}

  Let $(x,y,x',y')\in\Diam_2$ be a face of $\Diam$, the Cauchy-Riemann
  equation for $f\cdot\alpha$, on the couple $(x,x')$ and $(y,y')$ is
  the nullity of:
\begin{align}
  \hidewidth \frac{\int\limits_{(y,y')}f\cdot\alpha}{\ell(y,y')}-
  &i\frac{\int\limits_{(x,x')}f\cdot\alpha}{\ell(x,x')}\notag\\ 
  =&\frac1{\ell(y,y')}
  \bigl(\frac{f(x)+f(y)}2\int\limits_{\hidewidth(y,x)\hidewidth}\alpha+
  \frac{f(x)+f(y')}2\int\limits_{\hidewidth(x,y')\hidewidth}\alpha
  +\frac{f(x')+f(y)}2\int\limits_{\hidewidth(y,x')\hidewidth}\alpha+
  \frac{f(x')+f(y')}2\int\limits_{\hidewidth(x',y')\hidewidth}\alpha
  \bigr)\notag\\ 
  \hidewidth-i\frac1{\ell(x,x')}& 
  \bigl(\frac{f(x)+f(y)}2\int\limits_{\hidewidth(x,y)\hidewidth}\alpha
  +\frac{f(x')+f(y)}2\int\limits_{\hidewidth(y,x')\hidewidth}\alpha
  +\frac{f(x)+f(y')}2\int\limits_{\hidewidth(x,y')\hidewidth}\alpha
  +\frac{f(x')+f(y')}2\int\limits_{\hidewidth(y',x')\hidewidth}\alpha
  \bigr)\notag\\ 
  =&\bigl(\int\limits_{\hidewidth(y,x)\hidewidth}\alpha+
  \int\limits_{\hidewidth(y',x')\hidewidth}\alpha\bigr)
  \frac{f(y')-f(y)}{\ell(y,y')},\notag
\end{align}
after having developed, used the holomorphy of $\alpha$, then the
holomorphy of $f$.$\Box$

So to be able to construct out of the holomorphic $1$-forms $dZ$ given
by local flat isometries, and a holomorphic function a holomorphic 
$1$-form
$fdZ$, we have to impose that for each face $(x,y,x',y')\in\Diam_2$,
$Z(x)-Z(y)=Z(y')-Z(x')$. Geometrically, it means that each face of the
graph $\Diam$ is mapped by $Z$ to a parallelogram in $\mathbb{C}$. But
as the diagonals of this parallelogram are orthogonal, it is a
lozenge (or rhombus, or diamond).

\begin{definition}
  A double $(\Lambda,\ell)$ of a Riemann surface $\Sigma$ is {\bf
    critical} if it is semi-critical and each face of $\Diam_2$ are
  lozenges. Let $\delta(\Lambda)$ be the common length of their sides.
\end{definition}

\begin{remark} \rm
  This has an intrinsic meaning on $\Sigma$, the faces of $\Diam$ are
  genuine lozenges on the surface and every edge of $\Lambda$ can be
  realized by segments of length given by $\ell$, two dual edges being
  orthogonal segments.

  Another equivalent way to look at criticality can be useful: a
  double $(\Lambda,\ell)$ is critical if there exists an application
  $Z:\widetilde{\Sigma\setminus P}\to\mathbb{C}$ from the universal
  covering of the punctured surface $\Sigma\setminus P$ for a finite
  set $P\subset\Lambda_0$ into $\mathbb{C}$ identified to the oriented
  Euclidean plane $\mathbb{R}^2$ such that
\begin{itemize}
\item the image of an edge $a\in\tilde\Lambda_1$ is a linear segment
  of length $\ell(a)$,
\item two dual edges are mapped to a direct orthogonal basis,
\item $Z$ is an embedding out of the vertices,
\item there exists a representation $\rho$ of the fundamental group
  $\pi_1(\Sigma\setminus P)$ into the group of isometries of the plane
  respecting orientation such that,
  $$\forall\gamma\in\pi_1(\Sigma\setminus P),
  Z\circ\gamma=\rho(\gamma)\circ Z,$$
\item and the lengths of all the segments corresponding to the edges
  of $\Diam$ are all equal to the same $\delta>0$.
\end{itemize}



The criticality of a double map gives a coherent system of angles
$\phi$ in $(0,\pi)$ on the unoriented edges of $\Lambda$, $\phi(x,x')$
is the angle in the lozenge for which $(x,x')$ is a diagonal, at the
vertex $x$ (or $x'$). They verify that for every lozenge, the sum of
the angles on the dual diagonals is $\pi$. Then the conic angle at a
vertex is given by the sum of the angles over the incident edges.
\end{remark}   

Every discrete conformal structure $(\Lambda,\ell)$ defines a
conformal structure on the associated topological surface by pasting
lozenges together according to the combinatorial data (though most of 
the 
vertices will be conic singularities). Conversely,
 \begin{theorem}  \label{th:crit}
  Every closed Riemann surface accepts a critical map.
\end{theorem}
\begin{demo}{\ref{th:crit}}
  We first produce critical maps for cylinders of any modulus:
  Consider a row of $n$ squares and glue back its ends to obtain a
  cylinder, its modulus, the ratio of the square of the distance from
  top to bottom by its area is $\frac 1n$.

  Stacking $m$ such rows upon each other, one gets a cylinder of
  modulus $\frac mn$.

  Squares can be bent   into lozenges
yielding a
  continuous family of cylinders of moduli ranging from zero to
  $\frac 2n$ (see Fig.~\ref{fig:panto}). Hence we can get cylinders
  of any modulus.
\begin{figure}[htbp]
\begin{center}\input{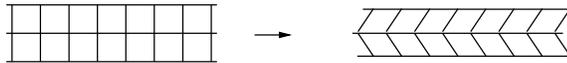}
\end{center}
\caption{Two bent rows.}         \label{fig:panto}

\end{figure}

Dehn twists can be performed on these critical cylinders, see
Fig.~\ref{fig:Dehn}.
\begin{figure}[htbp]
\begin{center}\input{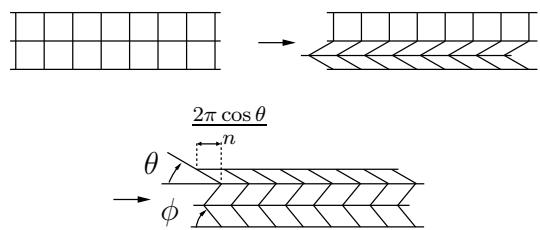}
\end{center}
\caption{Performing a Dehn twist.}         \label{fig:Dehn}
\end{figure}

Gluing three cylinders together along their bottom ($n$ has to be
even), one can produce trinions of any modulus (see
Fig.~\ref{fig:trinion}) and these trinions can be glued together
according to any angle. Hence, every Riemann surface can be so
produced~\cite{Bus}.
\begin{figure}[htbp]
\begin{center}\input{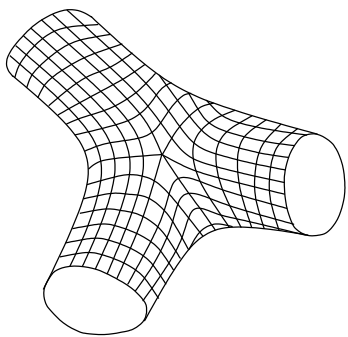}
\end{center}
\caption{Gluing three cylinders into a trinion.}         
\label{fig:trinion}
\end{figure}
\end{demo}

\begin{remark}
An equilateral surface is a Riemann surface which can be triangulated
by equilateral triangles with respect to a flat metric with conic
singularities. Equilateral surfaces are the algebraic curves over
$\bar{\mathbb{Q}}$~\cite{VoSh} so are dense among the Riemann surfaces.
Cutting every equilateral triangle into nine, three times smaller,
triangles (see Fig.~\ref{fig:equilatri}), one can couple these
triangles by pairs so that they form lozenges, hence a critical map.
\begin{figure}[htbp]
\begin{center}\input{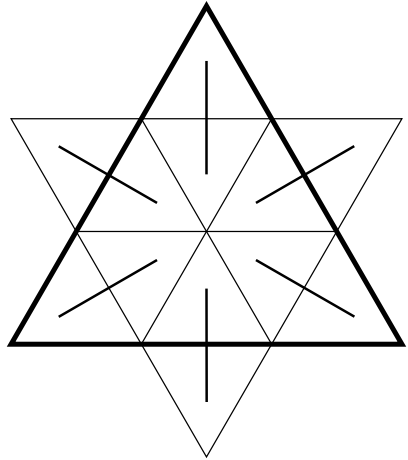}
\end{center}
\caption{An equilateral triangle cut in nine yielding lozenges.}
\label{fig:equilatri}

\end{figure}
\end{remark}

\begin{figure}[htbp]
\begin{center}\input{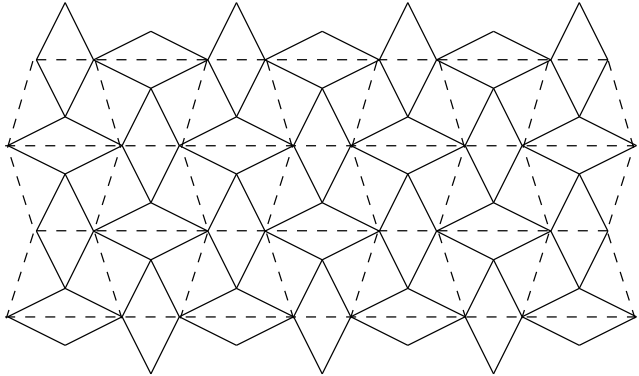}
\end{center}
\caption{A $1$-parameter family of critical deformations of the square
  lattice.}
\label{fig:carrecrit}

\end{figure}

\begin{figure}[htbp]
\begin{center}\input{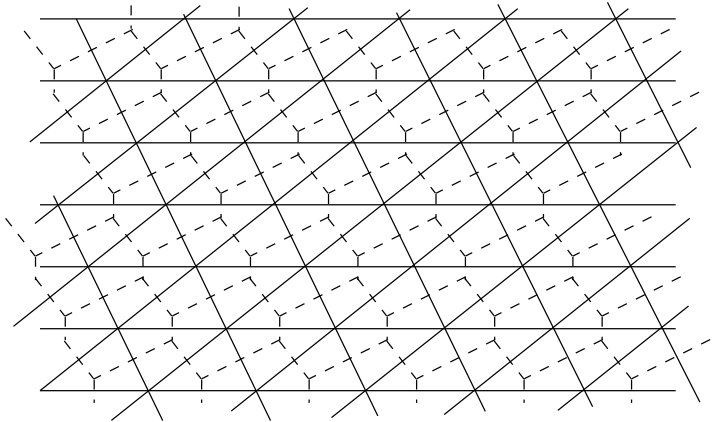}
\end{center}
\caption{A $2$-parameters family of critical deformations of the
  triangular/hexagonal lattices. This family, key to the solution of
  the triangular Ising model, induced Baxter to set up the Yang-Baxter
  equation~\cite{Bax}. Our notion of criticality fits beautifully into
  this framework.}
\label{fig:triangles}

\end{figure}

\begin{figure}[htbp]
\begin{center}\input{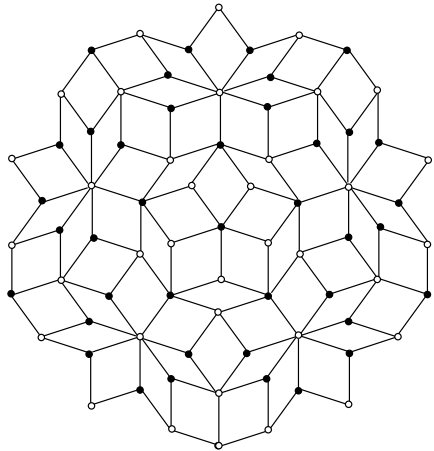}
\end{center}
\caption{The order $5$ Penrose quasi crystal.}
\label{fig:quasi}

\end{figure}

\begin{figure}[htbp]
\begin{center}\input{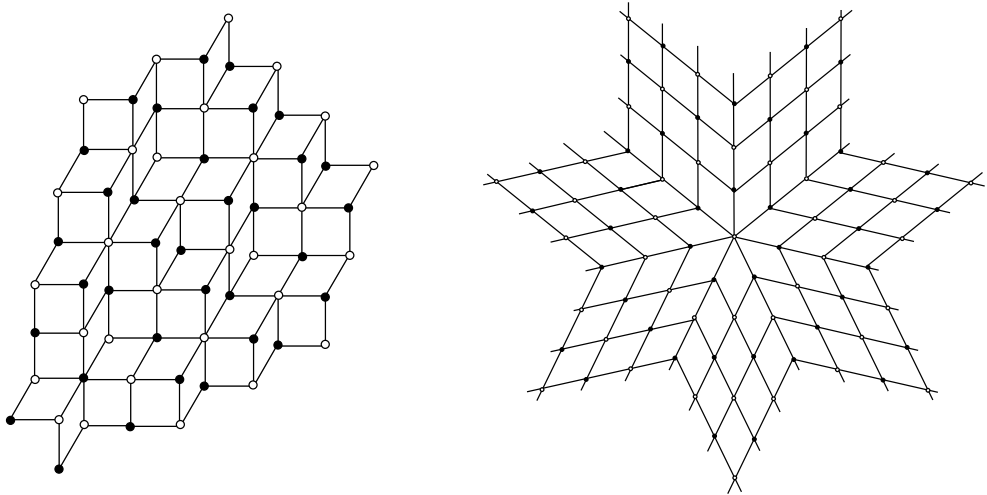}
\end{center}
\caption{Lozenge patchworks.}
\label{fig:patch}

\end{figure}

\begin{figure}[htbp]
\begin{center}\input{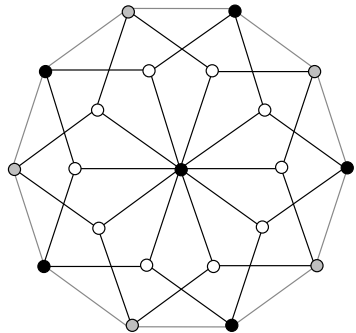}
\end{center}
\caption{Higher genus critical handlebody.}
\label{fig:losg2}

\end{figure}

\begin{figure}[htbp]
\begin{center}\input{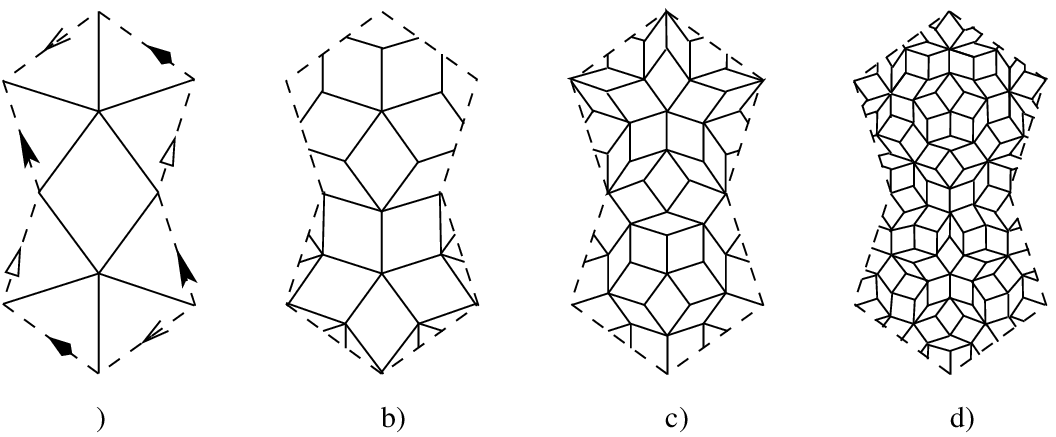}
\end{center}
\caption{Sequence of critical maps of a genus two handlebody using
Penrose inflation rule.}
\label{fig:pentag2}

\end{figure}

In Figures~\ref{fig:carrecrit}--\ref{fig:patch} are some examples of
critical decompositions of the plane. In Fig.~\ref{fig:losg2}, a
higher genus example, found in Coxeter~\cite{Cox1}, of the cellular
decomposition of a collection of handlebodies (the genus depends on
how the sides are glued pairwise) by ten regular pentagons, the centre
is a branched point of order three;  together with its dual, they form
a critical map. It is the case for any cellular
decomposition by just one regular tile when its vertices are
co-cyclic.
This decomposition gives rise to a critical sequence using the Penrose
inflation rule~\cite{GS87}. Fig.~\ref{fig:pentag2} illustrates this
inflation rule sequence on a simpler genus two example where each
outer side has to be glued with the other parallel side.

\subsection{Physical interpretation}

 \begin{theorem} \label{th:statcrit} A translationally invariant
 discrete conformal structure $(\Lambda,\rho)$ on $\Lambda$ the double
 square or triangular/hexagonal lattices decomposition of the plane or
 the genus one torus, is critical and flat if and only if the Ising
 model defined by the interaction constants
 $K_e:=\frac12\mathrm{Arcsinh}\rho_e$ on each edge $e\in\Lambda_1$ is
 critical as usually defined in statistical mechanics~\cite{McCW}.
\end{theorem}
\begin{demo}{}
We prove it by solving another problem which contains these two
particular cases, namely the translationally invariant square lattice
with period two~\cite{Yam}. At a particular vertex, the flat critical
condition on the four conformal parameters is:
$$\sum_{i=1}^4\arctan\rho_i=\pi,$$ which is obviously invariant by 
all the
symmetries of the problem, including duality. When
$\rho_i=\rho_{i+2}$, we get the usual period one Ising model
criticality on the square lattice $$\sinh 2K_h\sinh 2K_v=1,$$ and
likewise when one of the four parameters degenerates to zero or
infinity, the three remaining coefficients fulfill $$\sinh 2K_I\sinh
2K_{II}\sinh 2K_{III}=
\sinh 2K_I+\sinh 2K_{II}+\sinh 2K_{III}$$ which is (a form of) the
criticality condition for the triangular/hexagonal Ising model. The 
case shown
in Fig.~\ref{fig:carrecrit} occurs when $\rho_1=\rho_3=1$, implying
$\rho_2\rho_4=1$.
\end{demo}

We see here that flat criticality, when the angles at conic
singularitites are multiples of $2\pi$, is more meaningful than
criticality in general.  This theorem is important because it
shows that statistical criticality is meaningful even at the finite
size level. It is well known~\cite{KW} that for lattices, it
corresponds to self-duality, which has a meaning for finite systems;
here we see that self-duality corresponds to a compatibility with
holomorphy. In a sense, our notion of criticality defines self-duality
for more complex graphs than lattices. Furthermore, we will see in
Sect.~\ref{sec:spin} that criticality implies the existence of a
discrete massless Dirac spinor, which is the core of the Ising
model. Although we saw that criticality implies a continuous limit
theorem, the thermodynamic limit is not necessary for criticality to
be detected, and to have an interesting meaning.

It is easy to produce higher genus flat critical maps and compute
their critical temperature, the examples in
Figures~\ref{fig:losg2}-\ref{fig:pentag2} have four kinds of
interactions corresponding to the diagonals of the two kinds of
quadrilateral tiles. They are critical when the angles of the
quadrilaterals are $\frac\pi 5,\frac{2\pi} 5,\frac{3\pi} 5$, and
$\frac{4\pi} 5$, corresponding to Ising interactions
\begin{equation}
\sinh 2K_n=\tan \frac{n\pi} {10}.
\end{equation}
The author had made no attempt to verify these values numerically.

A general way is, considering a critical genus one torus made up of a
translationally invariant lattice, to cut two parallel segments of
equal length and seam them back, interchanging their sides. This
creates two conic singularities where an extra curvature of $-2\pi$ is
concentrated at each point, yielding a genus two handlebody. Repeating
the process, we may produce critical handlebodies of arbitrarily large
genus if we start with a very fine mesh. One has to beware that our
continuous limit theorem applies only to fixed genus, it cannot grow
with the refinement of the mesh. This explains why the union-jack
lattice (the square lattice and its diagonals) or the three
dimensional Ising model, which can be modelled as a genus $mnp$
surface for a $2m\times 2n\times 2p$ cubic network, are beyond the
scope of our technique as far as a continuous limit theorem is
concerned. With this restriction in mind, we see that both the
existence and the value of a critical temperature is essentially a
local property and neither depends on the genus nor on the modulus of
the handlebody. It is not the case for more interesting quantities
such as the partition function, which can be obtained in principle
from the discrete Dirac spinor that criticality provides, defined in
Sect.~\ref{sec:spin}. But such a calculus is beyond the scope of
this article.

\begin{figure}[htbp]
\begin{center}\includegraphics{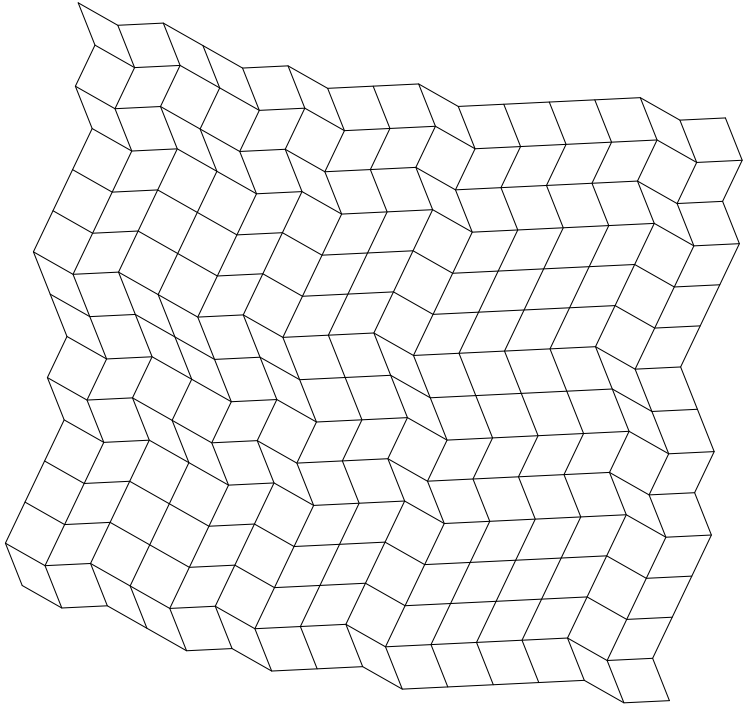}
\end{center}
\caption{The diamond graph of a critical labyrinth lattice.}
\label{fig:laby}
\end{figure}

\begin{figure}[htbp]
\begin{center}\input{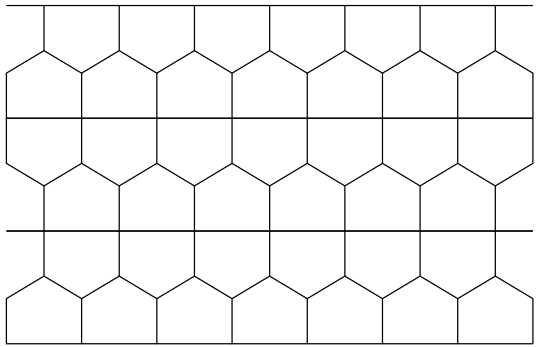}
\end{center}
\caption{The ``street'' lattice.}
\label{fig:street}
\end{figure}

Apart from the standard lattices, the critical temperature of other
well known graphs can be computed using our method, for example the
labyrinth~\cite{BGB}, whose diamond is pictured in
Fig.~\ref{fig:laby}, has the topology of the square lattice but has
five different interactions strengths controlled by two binary words,
labelling the columns and rows by 0's and 1's.
And also new ones such as the ``street graph'' depicted in
Fig.~\ref{fig:street}. Its double row transfer matrix appears to be
the product of three commuting transfer matrices, two triangular and
a square one.

  Other cases such as the
Kagom\'e~\cite{Syo} or more generally lattices of chequered
type~\cite{Uti} can be handled using a technique called electrical
moves~\cite{CdV96} which enables us to move around, and causes
appearing or disappearing conic singularities of a flat
metric. This will be the subject of a subsequent article, explaining
the relationship between discrete holomorphy, electrical moves and
knots and links.  These electrical moves act in the space of all the
graphs with discrete conformal structures in a similar way to that of
the Baxterisation processes in the spectral parameter space of an
integrable model (see~\cite{AdABM}). We are going to see that the link
with statistical mechanics is even deeper than simply pointing out a
submanifold of critical systems inside the huge space of all Ising
models, as the similarity with the continuous case extends to the
existence of a discrete Dirac spinor near criticality.

\subsection{Polynomial ring}
\begin{definition}
  Let $(\Lambda,\ell)$ be a critical map. In a given flat map $Z:U\to
  \mathbb{C}$ on the simply connected $U$, choose a vertex
  $z_0\in\Lambda_0$, and for a holomorphic function $f$, define the
  holomorphic functions $f^\dag$ and $f'$ by the following formulae:
  $$f^\dag(z):=\varepsilon(z)\bar f(z),$$ where $\bar f$ denotes the
  complex conjugate and $\varepsilon(\Gamma)=+1,
  \varepsilon(\Gamma^*)=-1$,
  $$f'(z):=\frac4{\delta^2}\left( \int_{z_0}^z f^\dag
    dZ\right)^\dag.$$
\end{definition}
See~\cite{Duf} for similar definitions. Notice that $f'$ is defined up
to $\varepsilon$ if one changes the base point.

\begin{proposition}
  Let $(\Lambda,\ell)$be  a critical map. In a given flat map $Z:U\to
  \mathbb{C}$ on the simply connected $U$, for every holomorphic
  function $f\in\Omega(\Lambda)$, $df=f'dZ$. We hence call $f'$ the
  derivative of $f$.
\end{proposition}
Consider an edge $(x,y)\in\Diam_1$, $x\in\Gamma_0, y\in\Gamma^*_0$,
\begin{align}
  f'(y)&=\frac4{\delta^2}\left( \int_{z_0}^x f^\dag dZ+\int_{x}^y
    f^\dag dZ\right)^{\dag_y}\notag\\ &=-f'(x)+
  \frac4{\delta^2}\left(\frac{\bar f(x)-\bar
      f(y)}2(Z(y)-Z(x))\right)^{\dag_y}\notag\\ 
  &=-f'(x)-\frac2{\delta^2} (f(x)-f(y))(\bar Z(y)-\bar Z(x)).\notag
\end{align}
So $\int_{(x,y)}f'dZ=- \frac{f(x)-f(y)}{\delta^2}(\bar Z(y)-\bar
Z(x))(Z(y)- Z(x))=f(y)-f(x)$.$\Box$

\begin{definition} Let $U$ be a simply connected flat region and  
$z_0\in U$.
  Define inductively the holomorphic functions $Z^k(z):=\int_{z_0}^z
  \frac1kZ^{k-1}dZ$ given $Z^0:=1$.  As the space of holomorphic
  functions on $U$ is finite dimensional, these functions are not
  free; let $P_U$ be the minimal polynomial such that
  $P_U(Z)=Z^n+\ldots=0$. 
 \end{definition}


 
\begin{conj}
  The space of holomorphic functions on $U$, convex, is isomorphic to
  $\mathbb{C}[Z]/P_U$.
\end{conj}

We won't define here the notion of convexity, see~\cite{CdV96}.
The question is whether the set $(Z^k)$ generates the space of
holomorphic functions.  The problem is that zeros are not localised,
and as the power of $Z^k$ increases, the set of its zeros spread on
the plane and get out of $U$. Figure
~\ref{fig:grdg1516} is an example on the unit square lattice with $U$ 
the square
$[-10,10]\oplus [-10,10]i$, the degree increases with $k$ until $16$
where four zeros get out of the square. So a definition of the degree
of a function by a Gauss formula is delicate.
\begin{figure}[htbp]
\begin{center}
  \input{}
\end{center}
\caption{The zeros of $Z^{16}$ get out of the square $[-10,10]\oplus
  [-10,10]i$.} \label{fig:grdg1516}

\end{figure}

\section{Dirac Equation}                        \label{sec:spin}

Although we believe our theory can be applied to a lot of different
problems, our motivation was to shed new light on statistical
mechanics and the Ising model in particular. This statistical model
has been linked with Dirac spinors since the work of Kaufman~\cite{K}
and Onsager and Kaufman~\cite{KO}. We refer among others
to~\cite{McCW81,SMJ,KC}. Hence we are interested in setting up a Dirac
equation in the context of discrete holomorphy. To achieve this goal
we first have to define the discrete analogue of the fibre bundle on
which spinors live. We therefore have to define a discrete spin
structure. Physics provides us with a geometric definition~\cite{KC}
based on paths in a certain $\mathbb{Z}_2$-homology, that we
generalise to our need (higher genus, boundary, arbitrary
topology). We begin by showing that such an object in the continuum is
indeed a spin structure, then define the discrete object. We then set
up the Dirac equation for discrete spinors, show that it implies
holomorphy and that the existence of a solution is equivalent to
criticality. The Ising model gives us an object which satisfies the
discrete Dirac equation, namely the fermion, $\Psi=\sigma\mu$ as
defined in~\cite{KC}, corresponding to a similar object defined
previously by Kaufman~\cite{K}. It fulfills the Dirac equation at
criticality, but also off criticality, corresponding to
a \textbf{massive} Dirac spinor. We will end this article by 
describing
off-criticality, as defined by the author's Ph.D. advisor, Daniel
{Bennequin}.

\subsection{Universal spin structure}

A spin structure~\cite{Mil} on a principal fibre bundle $(E,B)$ over a
manifold $B$, with $\mathrm{SO}(n)$ as a structural group, is a
principal fibre bundle $(E',B)$, of structural group
$\mathrm{Spin}(n)$, and a map $f:E'\to E$ such that the following
diagram is commutative:
$$
\begin{array}{cclcclc}
  E'&\times&\mathrm{Spin}(n)&\to&E'\\ &&&&&\searrow\\ 
  &\downarrow&f\times\lambda&&\downarrow f&&B\\ &&&&&\nearrow\\ 
  E&\times&\mathrm{SO}(n)&\to&E
\end{array}
$$ where $\lambda$ is the standard $2$-fold covering homomorphism from
$\mathrm{Spin}(n)$ to $\mathrm{SO}(n)$.

In this paper we consider only spin structures on the tangent bundle
of a surface.  On a generic Riemann surface $\Sigma$, there is not a
canonical spin structure.  We are going to describe a surface
$\hat\Sigma$, $2^{2-\chi(\Sigma)}$-fold covering of $\Sigma$, on which
there exists a preferred spin structure. It allows us to define every
spin structure on $\Sigma$ as a quotient of this universal spin
structure. We will treat the continuous case and then the discrete
case.

\begin{definition} Let $\Sigma$ be a differentiable surface with a 
base-point $y^0$; 
$\hat\Sigma$ is the set of pairs $(z,[\lambda]_2)$, where $z\in\Sigma$
is a point and $[\lambda]_2$ the homology of a path $\lambda$ from
$y^0$ to $z$ considered in the relative homology
$H_1(\Sigma,\{y^0,z\})\otimes\mathbb{Z}_2$.
\end{definition}

$\hat\Sigma$ is the $2^{2-\chi(\Sigma)}$ covering associated to the
intersection $H$ of the kernels of all the homomorphisms from
$\pi_1(\Sigma)$ to $\mathbb{Z}_2$, that is to say the quotient of the
universal covering by the subgroup $H\subset \pi_1(\Sigma)$ of loops
whose homology is null modulo two.

\begin{figure}[htbp]
\begin{center}
  \input{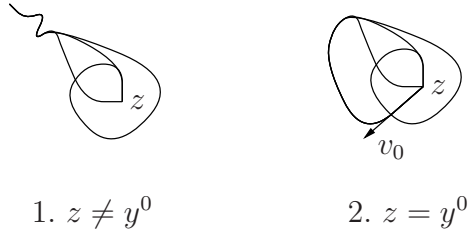}
\end{center}
\caption{Paths of different classes with
  respect to $\sim_z$ for $z\not = y^0$ and $z= y^0$.}
\label{fig:simz}

\end{figure}

 Choose $v_0$ a tangent vector at $y^0$. For each point $z\in\Sigma$,
 define $\Sigma_z:=\Sigma\setminus\{
 y^0,z\}\sqcup\mathbb{S}^1\sqcup\mathbb{S}^1$, the blown up of
 $\Sigma$ at $y^0$ and $z$ (add only one circle in the case
 $y^0=z$). Consider the set of oriented paths in $\Sigma_z$, from the
 point corresponding to the vector $v_0$ at $y^0$ to the directions at
 $z$ (the vector $v_0$ is needed only when $z=y^0$). Define an
 equivalence relation $\sim_z$ (see Fig.~\ref{fig:simz}) on this set
 by stating that two paths $\lambda$, $\lambda'$ are equivalent if and
 only if $\lambda-\lambda'$ is a cycle and $[\lambda-\lambda']_2=0$ in
 the homology $H_1(\Sigma\setminus \{z\},\mathbb{Z}_2)$.

\begin{definition}
The \textbf{universal spin structure} $\mathcal{S}$ of $\Sigma$ is the
set of pairs $(z,[\lambda]_{\sim_z})$, with $z\in\Sigma$ and
$[\lambda]_{\sim_z}$ the $\sim_z$-equivalence class of the path
$\lambda$ from $y^0$ to $z$ in $\Sigma_z$.
\end{definition}

\begin{theorem}                    \label{prop:S2S1}
  $\mathcal{S}$ is a spin structure on $\hat\Sigma$ and is the only
  one such that the action of the fundamental group $\pi_1(\Sigma)$ on
  $\hat\Sigma$ can be lifted to. Moreover it is the pull-back of any
  spin structure on $\Sigma$.
\end{theorem}
\begin{demo}{}
  The proof is in three steps, we check that $\mathcal{S}$ is a spin
  structure, we define a spin structure $\mathcal{S}_0$ through group
  theory and we show that both are equal to a third spin structure
  $\mathcal{S}_1$.

  There is an obvious projection from $\mathcal{S}$ to $\hat\Sigma$
  defined by $(z,[\lambda]_{\sim_z})\mapsto(z,[\lambda]_2)$. The fibre
  of this projection at $(z,[\lambda]_2)$ is the set of
  $\sim_z$-equivalence classes of paths from $y^0$ to the blown-up
  circle at $z$. To each class is associated the tangent direction at
  $z$ so $\mathcal{S}_z$ is a covering of $ST_z(\hat\Sigma)$. As
  $H_1(\Sigma\setminus\{z\},\mathbb{Z}_2)$ is $2^{3-\chi(\Sigma)}$
  dimensional (a loop around $z$ is not homologically trivial), for
  each point in $ST(\hat\Sigma)$, there are two different lifts. The
  path in $\mathcal{S}_z$ corresponding to turning around $z$ once
  yields the $\mathbb{Z}_2$-deck transformation. Hence $\mathcal{S}$
  is a spin structure on $\hat\Sigma$. 

  Let $G:=\pi_1(\Sigma)$ and $G':=\pi_1(ST\Sigma)$; the
  $\mathbb{S}^1$-fibre bundle $ST(\Sigma)\to \Sigma$ induces a short
  exact sequence $\mathbb{Z}\inj G'\surj G$. Every double covering of
  $ST\Sigma$ is defined by the kernel $S'$ of an homomorphism $u$
  from $G'$ to $\mathbb{Z}/2$, moreover, for $S'$ to be a spin
  structure, its intersection with the subgroup $\mathbb{Z}$ must be
  $2\mathbb{Z}$.

  Likewise, the fibration $\hat\Sigma\to \Sigma$ implies that the
  fundamental group $H':=\pi_1(ST\hat\Sigma)$ of the directions bundle
  of $\hat\Sigma$ is the subgroup of $G'$ over
  $H:=\pi_1(\hat\Sigma)$,

\begin{equation}
\begin{array}{ccccc}
\mathbb{Z}&\to&H'&\to&H\\
\downarrow&&\downarrow&&\downarrow\\
\mathbb{Z}&\to&G'&\to&G
\end{array},
\end{equation}

The intersection of the subgroups $H'$ and $S'$ is a well defined spin
structure $\mathcal{S}_0$ on $\hat\Sigma$: Indeed, consider another
spin structure $S''=\mathrm{Ker~}(v:G'\to\mathbb{Z}/2)$ on $\Sigma$,
its intersection with $\mathbb{Z}$ is $2\mathbb{Z}$ hence the kernel
of $u-v$ contains the whole subgroup $\mathbb{Z}$, that is to say
$u-v$ comes from a homomorphism of $G$ to $\mathbb{Z}/2$ and we have
$S''\cap H'=S'\cap H'$.  In other words, $\mathcal{S}_0$ is the unique
spin structure on $\hat\Sigma$ which is the pull-back of a spin
structure on $\Sigma$ and it is the pull-back of any spin structure.

Let $z\in\Sigma$ be a point, consider the set of paths in $ST\Sigma$ 
from
the base point $(y^0,v^0)$ to any direction at $z$. Consider on this
set the equivalence relation $\sim'_z$ defined by fixed extremities
$\mathbb{Z}/2$-homology. The class $[\lambda]_{\sim'_z}$ of a path
$\lambda$ from $(y^0,v^0)$ to $(z,v)$ is its homology class in
$H_1(ST\Sigma,\{ (y^0,v^0),(z,v)\})\otimes\mathbb{Z}/2$. The
projection $ST\Sigma\surj\Sigma$ splits $H_1(ST\Sigma,\{
(y^0,v^0),(z,v)\})\otimes\mathbb{Z}/2$ into
\begin{equation}
  \label{eq:split}
 \mathbb{Z}/2\to H_1(ST\Sigma,\{ (y^0,v^0),(z,v)\})\otimes\mathbb{Z}/2
\to H_1(\Sigma,\{ y^0,z\})\otimes\mathbb{Z}/2,
\end{equation}
hence the set $\mathcal{S}_1$ of pairs $(z,[\lambda]_{\sim'_z})$ for
all points $z\in\Sigma$ and all paths $\lambda$, is a spin structure
on $\hat\Sigma$.

Let $S'$ be a spin structure on $\Sigma$, it defines an element in
$\mathbb{Z}/2$ for each loop in $ST\Sigma$. So each path in $ST\Sigma$
beginning at $(y^0,v^0)$ defines, through the
splitting~\ref{eq:split}, an element in $\mathcal{S}_1$ which is then
the pull-back of $S'$ to $\hat\Sigma$, hence
$\mathcal{S}_0=\mathcal{S}_1$.

On the other hand $\mathcal{S}=\mathcal{S}_1$ because there is a
continuous projection from $\mathcal{S}$ to $\mathcal{S}_1$: For an
element $(z,[\lambda]_{\sim_z})$, consider a $C^1$-path $\lambda\in
\Sigma$ representing the class. Lift it to a path in $ST\Sigma$ by the
tangent direction at each point, its class $[\lambda]_{\sim'_z}$ only
depends on $[\lambda]_{\sim_z}$ and gives us an element in
$\mathcal{S}_1$.
\end{demo}

\subsection{Discrete spin structure}

\begin{definition}
Let $\Upsilon$ be a cellular complex of dimension two, 
a \textbf{spin structure} on $\Upsilon$ is a {\em graph} $\Upsilon'$, 
double
cover of the $1$-skeleton of $\Upsilon$ such that the lift of the
boundary of  every
face is a non-trivial double cover. They are considered up to
isomorphisms.
Let $S_D$ be the set of such spin structures.

A \textbf{spinor} $\psi$ on $\Upsilon'$ is an equivariant complex
function on $\Upsilon'$ regarding the action of $\mathbb{Z}/2$, that
is to say, for all $\xi\in\Upsilon'_0$, $\psi(\bar\xi)=-\psi(\xi)$ if
$\bar\xi$ represents the other lift.
\end{definition}

\begin{remark} \rm
  Usually, a spinor field is a section of a spinor bundle, that is to
  say a square root of a tangent {\em vector} field. Here, we consider
  square roots of {\em covectors}; we should say {\em cospinors}.

  A discrete spin structure is encoded by a representation of the
  cycles of $\Upsilon$, $Z_1(\Upsilon):=\mathrm{Ker~}\partial\cap
  C_1(\Upsilon)$, into $\mathbb{Z}/2$ which associates to $\gamma\in
  Z_1(\Upsilon)$, the value $\mu(\gamma)=0$ if it can be lifted in
  $\Upsilon'$ to a cycle and $\mu(\gamma)=1$ if it can not. By
  construction, the value of the boundary of a face is $1$ and the
  value of a cycle which is the boundary of a $2$-chain of $\Upsilon$
  is the number of faces enclosed, modulo two.
\end{remark}

We are going to show that this structure is indeed a good notion of
discrete spin structure. First, there are as many discrete spin 
structures on
a surface as there are in the continuous case:
   
\begin{proposition}
  On a closed connected oriented genus $g$ surface $\Sigma$, the set
  $S_D$ of inequivalent discrete spin structures of a cellular
  decomposition $\Upsilon$ is of cardinal $2^{2g}$. The space of
  representations of the fundamental group of the surface into
  $\mathbb{Z}/2$ acts freely and transitively on $S_D$.
\end{proposition}
  We explicitly build discrete spin structures and count them:
  Let $T$ be a maximal tree of $\Upsilon$, that is to say a
  sub-complex of dimension one containing all the vertices of
  $\Upsilon$ and a maximal subset of its edges such that there is no
  cycle in $T$. 
Choose $2g$ edges $(e_k)_{1\leq k\leq 2g}$ in
  $\Upsilon\setminus T$ such that the $2g$ cycles $(\gamma_k)\in
  Z_1(\Upsilon)^{2g}$ extracted from $(T\cup e_k)_{1\leq k\leq 2g}$
  form a basis of the fundamental group of $\Sigma$ (and $\Upsilon$).
Let $T_+:=T\cup_k e_k$ and consider $T'$, the sub-complex of the dual
  $\Upsilon^*$ formed by all the edges in $\Upsilon^*$ not crossed by
  $T_+$. It is a maximal tree of $\Upsilon^*$. Likewise we define
  $T'_+:=T'\cup_k e^*_k$.

We construct inductively a spin structure $\Upsilon'$: its first 
elements are a
  double copy of $T$ and we add edges without any choice to make
  as we take leaves out of $T'_+$. When only cycles are left, a choice
  concerning an edge $e_k$ has to be taken, opening a cycle in
  $T'_+$. The process goes on until $T'_+$ is empty. 

  These choices are
completely encoded by a representation $\mu$ such as in the remark,
and the $2g$ values $(\mu(\gamma_k))_{1\leq k\leq 2g}$ determine the
spin structure.  On the other hand, this representation defines the
spin structure and there are $2^{2g}$ such different representations.
Hence the choices of the maximal tree and the edges $e_k$ are
irrelevant.
 
Because a cycle in $\Upsilon$ belongs to a class in the fundamental
group of the surface (up to a choice of a path to the base
point, irrelevant for our matter), the representations of the 
fundamental group
into $\mathbb{Z}/2$ obviously act on spin structures: A representation
$\rho:\pi_1(\Sigma)\to \mathbb{Z}/2$ associates to a spin structure
defined by a representation $\mu:Z_1(\Upsilon)\to \mathbb{Z}/2$, the
spin structure defined by the representation $\rho(\mu)$ such that
$\rho(\mu)(\gamma):=\mu(\gamma)+\rho([\gamma])$, where
$[\gamma]\in\pi_1(\Sigma)$ is the class of the cycle $\gamma$ in the
fundamental group. This action is clearly free, and transitive because
the set of representations is of cardinal $2^{2g}$.
$\Box$ 

Given $\Lambda=\Gamma\sqcup\Gamma^*$ a double cellular decomposition,
we introduce a cellular decomposition which is the discretised
version of the tangent directions bundle of both $\Gamma$ and 
$\Gamma^*$:

\begin{definition}
  The {\em triple graph} $\Upsilon$ is a cellular complex
  whose vertices are unoriented edges of $\Diam$,
  $\Upsilon_0=\left\{\{x,y\}/(x,y)\in\Diam_1\right\}$. Two vertices
  $\{x,y\},\{x',y'\}\in\Upsilon_0$ are neighbours in $\Upsilon$ iff 
the
  edges $(x,y)$ and $(x',y')$ are incident (that is to say $x=x'$ or
  $x=y'$ or $y=x'$ or $y=y'$), and they bound a common face of
  $\Diam$. There are two edges in $\Upsilon$ for each
  edge in $\Lambda$.  For this to be a cellular decomposition of the
  surface in the empty boundary case, one needs to add faces of
  three types, centred on vertices of $\Gamma$, of $\Gamma^*$ and on
  faces of $\Diam$ (see Fig.~\ref{fig:TripleFac}).
\begin{figure}[htbp]
\begin{center}
  \input{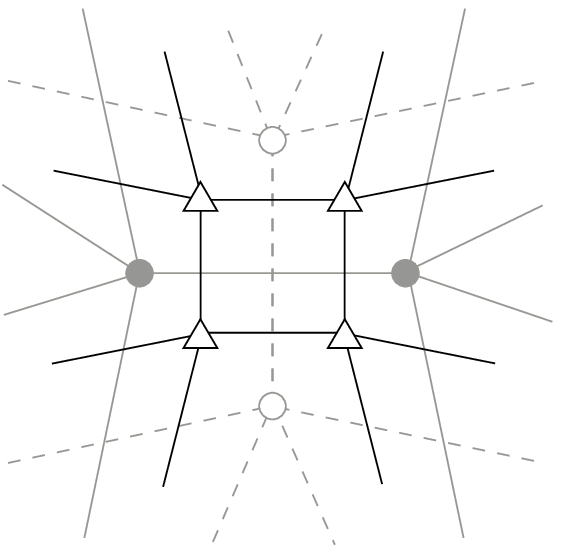}
\end{center}
\caption{The triple graph $\Upsilon$.}  \label{fig:TripleFac}

\end{figure}
\end{definition}

\begin{remark} \rm
  The topology of the usual tangent directions bundle is not at all
  mimicked by the incidence relations of $\Upsilon$, the former is 
$3$  dimensional and the latter is a $2$-cellular complex.
 \end{remark}

 Let
$(x^0,y^0)\in\Diam_1$ be a given edge. All the complexes $\Gamma,
\Gamma^*, \Diam, \Upsilon$ are lifted to $\hat\Sigma$.
  
\begin{definition}
  The \textbf{discrete universal spin structure} $\hat\Upsilon'$ is 
the
  following $1$-complex: Its vertices are of the form
  $((x,y),[\gamma^y_{y^0}])$, where $(x,y)\in\Upsilon_0$ is a pair of
  neighbours in $\Diam$ and $\gamma^y_{y^0}$ is a path from $y^0$ to
  $y$ on $\Gamma^*$, avoiding the faces $x^*$ and $x^{0*}$. We are
  interested only in its relative homology class modulo two, that is
  to say $[\gamma^y_{y_0}]\in H_1(\Gamma^*\setminus 
x^*,\{y^0,y\})\otimes\mathbb{Z}_2$.  We will denote a point by
  $((x,y),\gamma^y_{y^0})$ and identify it with
  $((x,y),{\gamma'}^y_{y^0})$ whenever $\gamma^y_{y^0}$ and
  ${\gamma'}^y_{y^0}$ are homologous. 

Two points
  $((x,y),\gamma^y_{y^0})$ and $((x',y'),{\gamma}_{y'}^{y^0})$ are
  neighbours in $\hat\Upsilon'$ if
  \begin{itemize}
  \item $x=x'$, $(y,y')\in\Gamma_1^*$ and
    ${\gamma}_{y}^{y^0}-{\gamma}_{y'}^{y^0}+(y,y')$ is homologous to
    zero in $H_1(\Gamma^*\setminus x^*)\otimes\mathbb{Z}_2$,
  \item $y=y'$, $(x,x')\in\Gamma_1^*$ and
    ${\gamma}_{y}^{y^0}-{\gamma}_{y'}^{y^0}$ is homologous to zero in
    $H_1(\Gamma^*\setminus x^*)\otimes\mathbb{Z}_2$.
  \end{itemize}
 \end{definition}

\begin{figure}[htbp]
\begin{center}
  \hfill \input{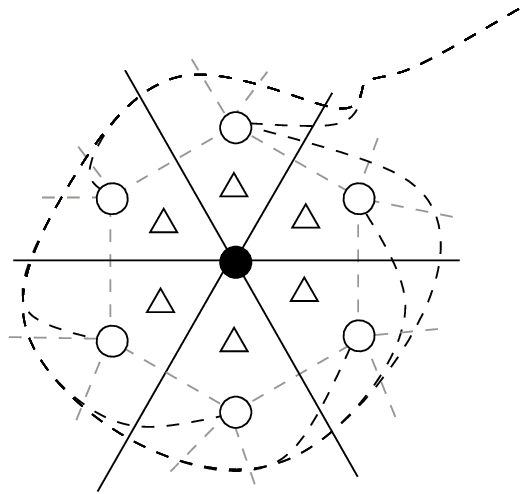} \hfill
  \input{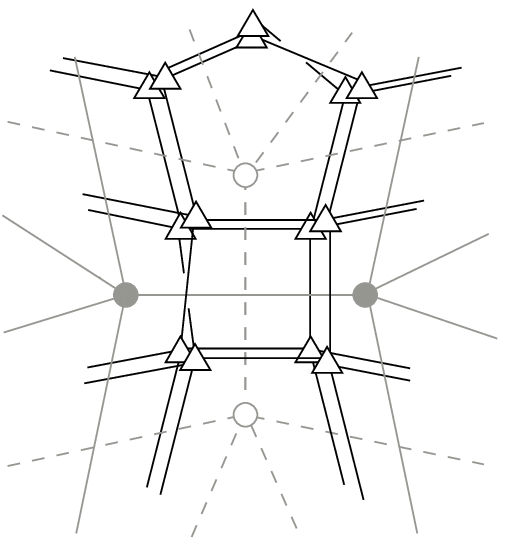}\hfill
\end{center}
\caption{Double covering around faces of $\Upsilon$.}   
\label{fig:TripleTilde}

\end{figure}

 $\hat\Upsilon'$ is a double covering of $\hat\Upsilon$ and it is
 connected around each face (see Fig.~\ref{fig:TripleTilde}). It is
 a discrete spin structure on $\hat\Upsilon$ in the sense defined
 above. Once a basis of the fundamental group $\pi_1(\Upsilon)$ is
 chosen, every representation of the homology group of $\Sigma$ into
 $\mathbb{Z}_2$ allows us to quotient this universal spin structure
 into a double covering of $\Upsilon$, yielding a usual spin structure
 $\Upsilon'$.

\subsection{Dirac equation}
A spinor changes sign between the two lifts in $\Upsilon'$ of a vertex
of $\Upsilon$, in other words it is multiplied by $-1$ when it turns
around a face. The faces of $\Upsilon$ which are centred on diamonds
are four sided.  We set up the \textbf{spin symmetry} equation for a
function $\zeta$ on $\Upsilon'_0$, on a positively oriented face
$(\xi_1,\xi_2,\xi_3,\xi_4)\in\Upsilon_2$ around a diamond, lifted
to an $8$-term cycle
$(\xi^+_1,\xi^+_2,\xi^+_3,\xi^+_4,\xi^-_1,\xi^-_2,\xi^-_3,\xi^-_4)\in
Z_1(\Upsilon')$:
\begin{equation}
  \label{eq:sym}
  \zeta(\xi^+_3)=i\zeta(\xi^+_1).
\end{equation}
It implies obviously that $\zeta$ is a spinor, that is to say
$\zeta(\xi^-_\bullet)=-\zeta(\xi^+_\bullet).$

The coherent system of angles $\phi$ given by a semi-critical
structure locally provides a spinor respecting the spin symmetry away
from conic singularities: Define half angles $\theta$ on oriented
edges of $\Upsilon$ in the following way: Each edge
$(\xi,\xi')\in\Upsilon_1$ cuts an edge $a\in\Lambda_1$, set
$\theta(\xi,\xi'):=\pm \frac{\phi(a)}2$ whether $(\xi,\xi')$ turns in
the positive or negative direction around the diamond. Choose a base
point $\xi_0\in\Upsilon'_0$, define $\zeta$ by $\zeta(\xi_0)=1$
and
\begin{equation}
  \label{eq:zetaxi}
  \zeta(\xi):=\exp i\sum_{\lambda\in\gamma} \theta(\lambda)
\end{equation}
for any path $\gamma$ from $\xi_0$ to $\xi$. The sum of the half
angles are equal to $\pi$ around the faces of $\Diam$ and half the
conic angle around a vertex, so if it is a regular flat point, we get
$\frac{2\pi}2=\pi$ again, hence $\zeta$ is a well defined spinor. As
diagonals of the faces of $\Diam$ are orthogonal, $\zeta$ fulfills the
spin symmetry. Moreover, if the conic angles are congruous to $2\pi$
modulo $4\pi$, $\zeta$ can be extended to any simply connected region;
if the fundamental group acts by translations, $\zeta$ is defined on
the whole $\Upsilon'$.

\begin{figure}[htbp]
\begin{center}
  \input{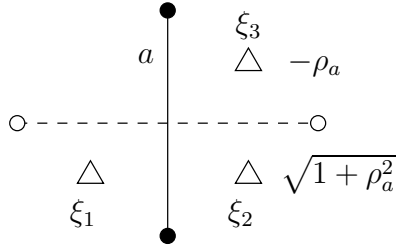}
\end{center}
\caption{The Dotsenko equation.}   \label{fig:DD}

\end{figure}

We are going to define a propagation equation which comes from
 the Ising model. It is fulfilled by the
fermion defined by Kaufman~\cite{K} which is known to converge to a
Dirac spinor near criticality. We will use the definition
$\psi=\sigma\mu$ given by Kadanoff and Ceva~\cite{KC}.  The Dirac
equation has a long history in the Ising model, beginning with the
work of Kaufman~\cite{K} and Onsager and Kaufman~\cite{KO}, we refer
among others to~\cite{McCW81,SMJ,KC}.  The equation that we need is
defined explicitly in~\cite{DD}, hence we will name it the {\bf
Dotsenko} equation, even though it might be found elsewhere in other
forms. It is fulfilled by the fermion at criticality as well as off
criticality. But this equation is only a part of the full Dirac
equation. For a function $\zeta$ on $\Upsilon'_0$, with the same
notations as before, and if $a\in\Lambda_1$ is the diagonal of the 
diamond,
between $(\zeta_2,\zeta_3)$ and $(\zeta_4,\zeta_1)$ (see
Fig.~\ref{fig:DD}):
\begin{equation}
  \label{eq:DD}  
\zeta(\xi_1^+)=\sqrt{1+\rho(a)^2}\zeta(\xi_2^+)-\rho(a)\zeta(\xi_3^+).
\end{equation}

A check around the diamond shows that it also implies that $\zeta$ is
a spinor: We write the Dotsenko equation in $\xi_2^+$ and $\xi_3^+$,
\begin{align}
\zeta(\xi_2^+)=&\sqrt{1+\rho(a^*)^2}\zeta(\xi_3^+)-
\rho(a^*)\zeta(\xi_4^+),
\notag\\ 
\zeta(\xi_3^+)=&\sqrt{1+\rho(a)^2}\zeta(\xi_4^+)-
\rho(a)\zeta(\xi_1^-),
\notag
\end{align}
hence, as $\sqrt{1+\rho(a)^2}\sqrt{1+\rho(a^*)^2}=\rho(a)+\rho(a^*)$,
\begin{align}
  \zeta(\xi_1^+)=&\rho(a^*)\zeta(\xi_3^+) -
  \sqrt{1+\rho(a^*)^2}\zeta(\xi_4^+)\notag\\ =&
  \rho(a^*)(\sqrt{1+\rho(a)^2}\zeta(\xi_4^+) -
  \rho(a)\zeta(\xi_1^-))-\sqrt{1+\rho(a^*)^2}\zeta(\xi_4^+)\notag\\ 
  =&-\zeta(\xi_1^-)\notag.
\end{align}

The \textbf{Dirac equation} is the conjunction of the
symmetry~(\ref{eq:sym}) and the Dotsenko~(\ref{eq:DD}) equations. We
will see that this same equation describes the massive and massless
Dirac equation, the mass measuring the distance from criticality.

Given two spinors $\zeta$, $\zeta'$, their pointwise product is no
longer a spinor but a regular function on $\Upsilon$. As there are two
edges in $\Upsilon$ for each edge in $\Lambda$, there is an obvious
averaging map from $1$-forms on $\Upsilon$ to $1$-forms on $\Lambda$:
We define $d_\Upsilon\zeta\zeta'\in C^1(\Lambda)$ by the following
formula, with the same notation as before,
$$2\int_{a} d_\Upsilon\zeta\zeta':= \zeta(\xi_3)
\zeta'(\xi_3)-\zeta(\xi_2) \zeta'(\xi_2) +\zeta(\xi_4)
\zeta'(\xi_4)-\zeta(\xi_1) \zeta'(\xi_1).
$$
$d_\Upsilon\zeta\zeta'$ is by definition an exact $1$-form on
$\Upsilon$ but its average {\em is not} {\it a priori} exact on
$\Lambda$.
\begin{figure}[htbp]
\begin{center}
  \input{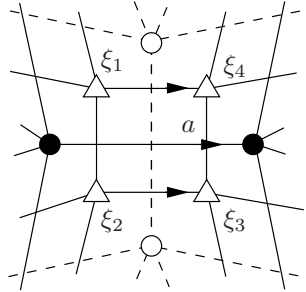}
\end{center}
\caption{The $1$-form on $\Lambda$ associated to two spinors.}   
\label{fig:zetazeta}
\end{figure}

\begin{proposition}
  If $\zeta$ and $\zeta'$ respect whether the spin symmetry or the
  Dotsenko equation, then $d_\Upsilon\zeta\zeta'$ is a closed
  $1$-form. If $\zeta$ is a Dirac spinor and $\zeta'$ fulfills the
  Dotsenko equation, then $d_\Upsilon\zeta\zeta'$ is holomorphic,
  $d_\Upsilon\bar\zeta\zeta'$ anti-holomorphic and every holomorphic
  $1$-form on $\Lambda$ can be written this way on a simply connected
  domain, uniquely up to a constant.
\end{proposition}

A sufficient condition for $d_\Upsilon\zeta\zeta'$ to be closed on 
$\Lambda$ is
that, with the same notations as above, $\zeta(\xi_3)
\zeta'(\xi_3)-\zeta(\xi_2) \zeta'(\xi_2)=\zeta(\xi_4)
\zeta'(\xi_4)-\zeta(\xi_1) \zeta'(\xi_1)$ because $\oint_{\partial
  y^*}d_\Upsilon\zeta\zeta'$ for a vertex $y\in\Lambda_0$ is a sum of
such differences on the edges of $\Upsilon$ around $y$.  This is so if
there exists a $2\times 2$-matrix $A$ such that
$$
\begin{pmatrix}
  \zeta(\xi_3^+)\\ \zeta(\xi_2^+)
\end{pmatrix}
= A
\begin{pmatrix}
  \zeta(\xi_4^+)\\ \zeta(\xi_1^+)
\end{pmatrix},
$$ a similar formula for $\zeta'$, and ${}^tA
\begin{pmatrix}
  1&0\\ 0 &-1
\end{pmatrix}
A=\begin{pmatrix} 1&0\\ 0 &-1
\end{pmatrix}
$. The solutions are of the form $A=
\begin{pmatrix}
  \epsilon\sqrt{1+\lambda^2}& \lambda\\ \epsilon \lambda&
  \sqrt{1+\lambda^2}
\end{pmatrix}$
for a complex number $\lambda\in\mathbb{C}$, $\epsilon=\pm 1$ and a
determination of $\sqrt{1+\lambda^2}$. This is the case for the spin
symmetry, $\lambda=-i, \epsilon=+1$ and for the Dotsenko equation,
$\lambda=\rho(a), \epsilon=-1, \sqrt{1+\lambda^2}>0$.

If $\zeta$ is a Dirac spinor and $\zeta'$ fulfills the Dotsenko 
equation,
then
\begin{align}
  \int_{a^*} d_\Upsilon\zeta\zeta'=& \zeta(\xi_4^+)\zeta'(\xi_4^+)
  -\zeta(\xi_3^+)\zeta'(\xi_3^+)\notag\\ 
=&i\zeta(\xi_2^+)(\sqrt{1+\rho(a)^2}\zeta'(\xi_3^+)-
\rho(a)\zeta'(\xi_2^+))  
-i\zeta(\xi_1^+)\zeta'(\xi_3^+)\notag\\  
=&i\zeta(\xi_2^+)(\sqrt{1+\rho(a)^2}\zeta'(\xi_3^+)-
\rho(a)\zeta'(\xi_2^+))
\notag\\  
&\;\;-i(\sqrt{1+\rho(a)^2}\zeta(\xi_2^+)-
\rho(a)\zeta(\xi_3^+))\zeta'(\xi_3^+)
\notag\\ =&i\rho(a)\left(\zeta(\xi_3^+)\zeta'(\xi_3^+)
    -\zeta(\xi_2^+)\zeta'(\xi_2^+)\right) =i\rho(a)\int_{a}
  d_\Upsilon\zeta\zeta'.\notag
\end{align}
So $d_\Upsilon\zeta\zeta'$ is holomorphic. Of course,
$d\bar\zeta\zeta'$ is anti-holomorphic. Conversely, if
$d_\Upsilon\zeta\zeta'$ is holomorphic with $\zeta$ a Dirac spinor,
then $\zeta'$ fulfills the Dotsenko equation.

Given a holomorphic $1$-form $\alpha\in\Omega^{(1,0)}(\Lambda)$,
define $\alpha_\Upsilon$ on $\Upsilon_1$ by the obvious map
$\int_{(\{x,y\},\{y,x'\})}\alpha_\Upsilon:=\int_{(x,x')}\alpha$. It is
a closed $1$-form on $\Upsilon$ because $\alpha$ is closed on
$\Lambda$, so there exists a function $a$ on any simply connected
domain of $\Upsilon_0$, unique up to an additive constant, such that
$d_\Upsilon a=\alpha_\Upsilon$. A check shows that the only spinors
$\zeta''$ such that $d_\Upsilon\zeta\zeta''=0$ on $\Lambda$ are the
one proportional to $\bar\zeta$. It is consistent with the fact that
the Dirac spinor is of constant modulus (see
Eq.~(\ref{eq:zetaxieiphi})). Hence the function $\zeta':=a/\zeta$
on $\Upsilon'$ is the unique spinor (up to a constant times
$1/\zeta\sim\bar\zeta$) such that
$d_\Upsilon\zeta\zeta'=\alpha$.$\Box$

Notice that for $\zeta$ a Dirac spinor, the holomorphic $1$-form
associated to it on $\Lambda$ is locally, for a given flat coordinate
$Z$, $d_\Upsilon\zeta\zeta=\lambda dZ$, with $\lambda\in\mathbb{C}$ a
certain constant.

\subsection{Existence of a Dirac spinor}
\begin{theorem}\label{th:Dirac}
  There exists a Dirac spinor on a double map iff it is critical for a
  given flat metric with conic angles congruous to $2\pi$ modulo
  $4\pi$ and such that the fundamental group acts by translations. The
  Dirac spinor is unique up to a multiplicative constant.
\end{theorem}
\begin{demo}{\ref{th:Dirac}}
  Let $\zeta$ be a non-zero Dirac spinor. Consider a positively 
oriented
  face $(\xi_1,\xi_2,\xi_3,\xi_4)\in\Upsilon_2$ around a diamond
  with diagonals $a,a^*$  as in Fig.~\ref{fig:DD}, lifted to an
  $8$-term cycle
  
$(\xi^+_1,\xi^+_2,\xi^+_3,\xi^+_4,\xi^-_1,\xi^-_2,\xi^-_3,\xi^-_4)\in
  C_1(\Upsilon')$. The equation
  $$e^{i\frac{\phi(a)}2}=\frac{\rho(a^*)+i}{\sqrt{1+\rho(a^*)^2}}$$
  defines an angle $\phi(a)\in(0,\pi)$ for every edge $a\in\Lambda_1$.

  The Dotsenko and symmetry equations combine into
\begin{equation}
  \label{eq:zetaxieiphi}
  \zeta(\xi_2^+)=\frac{\rho(a)+i}{\sqrt{1+\rho(a)^2}}\zeta(\xi_1^+).
\end{equation}
The fact that $\zeta$ is a spinor implies that, summing the four
angles around the diamond, we get $e^{i(\phi(a)+\phi(a^*))}=-1$. As
each angle is less than $\pi$, their sum is equal to $\pi$. The same
consideration around a vertex $x\in\Lambda_0$, yields $\exp
i\sum_{(x,x')\in\Lambda_1} \frac{\phi(x,x')}2 =-1$. So $\phi$ is a
coherent system of angles and the map is critical with conic angles
congruous to $2\pi$ modulo $4\pi$.

Conversely, given $\phi$ a coherent system of angles with conic angles
congruous to $2\pi$ modulo $4\pi$, the preceding construction
described by Eq.~(\ref{eq:zetaxi}) gives the only Dirac spinor.
\end{demo}

In this case, $dZ$ is a well defined holomorphic $1$-form on the whole
surface.

\begin{corol} Let $(\Lambda,\rho)$ be
  a discrete conformal structure and $P$ a set of vertices, containing
  among others the vertices $v$ such that the sum $\sum_{e\sim v}
  \mathrm{Arctan} \rho(e)$, summed over all edges $e$ incident to $v$,
  is greater than $2\pi$. The discrete conformal structure is critical
  with $P$ as conic singularities if and only if there exist Dirac
  spinors on every simply connected domain containing no point of $P$.
\end{corol}

We define in which sense a discrete spinor converges to
a continuous spinor. We don't define these spinors on specific spin
structures but rather on the universal spin structure $\mathcal{S}$.

Consider a sequence of finer and finer critical maps such as in
Theorem~\ref{th:lim}.  Choose a converging sequence of base points
$(x^0_k,y^0_k)\in{}^k\Upsilon_0$ on each critical map such that the
direction sequence
$(\frac{\overrightarrow{x^0_ky^0_k}}{d(x^0_k,y^0_k)})$ converges to a
tangent vector $(x^0,v^0)$.

Consider a sequence of points $(x_k,y_k)\in{}^k\Upsilon_0$, defining a
sequence of points $(x_k)$ converging to $x$ in $\Sigma$ and a
converging sequence of directions
$v=\mathrm{lim}\frac{\overrightarrow{x_ky_k}}{d(x_k,y_k)}$.  By
compacity of the circle, there exist such sequences for every point
$x\in\Sigma$ and the criticality implies that it is in at least three
directions for flat points, separated by angles less than $\pi$.

The different limits allow us to identify, after a certain rank, the
relative homology groups $H_1({}^k\Gamma^{*}\setminus
x^*_k,\{y^0_k,y_k\})\otimes\mathbb{Z}_2$ with
$H_1(\Sigma_x,\{(x^0,v^0),(x,v)\})\otimes\mathbb{Z}_2$, the classes of
paths in the blown-up of $\Sigma$ at $x^0$ and $x$.

\begin{definition}
  We will say that a sequence $(\zeta_k)_{k\in\mathbb{N}}$ of spinors
  converges if and only if, for any converging sequences,
  $\left((x_k,y_k)\in{}^k\Upsilon_0\right)_{k\in\mathbb{N}}$ defining 
a
  limit tangent vector, and
  $\left([\lambda_k]\right)_{k\in\mathbb{N}}$ of classes of paths in
  ${}^k\Gamma^{*}$ from $y^0_k$ to $y$, avoiding the face $x_k^*$, the
  sequence of values $(\zeta_k(x_k,[\lambda_k]))$ converges.
\end{definition}

\begin{remark} \rm
  It defines a continuous limit spinor $\zeta$ by equivariance: Let
  $x\in\hat\Sigma$, the set $D_x$ of directions in which there exist
  converging sequences of discrete directions is by definition a
  closed set. Let $u,v$ two boundary directions of $D_x$ such that the
  entire arc $A$ of directions between them is not in $D_x$.
  Consider $[(x,[\lambda]_x),(x,[\lambda']_x)]\subset\mathcal{S}$ a
  lift of $A$. The circle $\mathbb{S}^1$ acts on the directions, hence
  on the $\sim_x$-classes, let
  $\psi\in (0,\pi)$
  the angle such that
  $(x,e^{i\psi}[\lambda]_x)=(x,[\lambda']_x)$. Define
  $$\zeta(x,e^{i\phi}[\lambda]_x):=e^{iv(\phi)} 
\zeta(x,[\lambda]_x),$$
  where $v(\phi)=\frac \phi\psi
  \frac{\zeta(x,[\lambda']_x)}{\zeta(x,[\lambda]_x)}$ ($v(\phi)=\frac
  \phi 2$ for Dirac spinors).
\end{remark}

\begin{theorem}
Given a sequence of critical maps  such
as in Theorem~\ref{th:lim} with Dirac spinors on all of them,
 they can be normed so that they converge to
the usual Dirac spinor on the Riemann surface.
\end{theorem}
In a local flat map $Z$, the square of the discrete Dirac spinor on
${}^k\hat\Upsilon'$ is (up to a multiplicative constant) the $1$-form
$dZ$ evaluated on the edges. Hence their sequence converges.

\subsection{Massive Dirac equation, discrete fusion algebra and 
conclusions}

For completeness and motivation, we describe below the situation
off-criticality where elliptic integrals come into play, and
investigate a form of the discrete fusion algebra in the Ising
model. This work was done by Daniel Bennequin and will be the subject
of a subsequent article.

A massive system in the continuous theory is no longer conformal. In
the same fashion, Daniel Bennequin defined a massive discrete system 
of
modulus $k$ as a discrete double graph $(\Lambda,\rho)$ such that, for
each pair $(a,a^*)$ of dual edges,
\begin{equation}
\rho(a)\rho(a^*)=\frac1k.
\end{equation}

The massless case corresponds to $k=1$. We showed that
criticality was equivalent to a coherent system of angles $\phi(a)$
such as shown in Fig.~\ref{fig:systangsemi},  defined by $\tan
\frac{\phi(a)}2=\rho(a)$, and adding up to $2\pi$ at
each vertex of the double, except at conic singularities. The Dirac
spinor was constructed using the half angles
$\frac{\phi(a)}2$. Similarly, for every edge, we define the massive
``half angle'' $u(a)$ as the elliptic integral
\begin{equation}
u(a):=\int_0^{\frac{\phi(a)}2}\frac{d\varphi}{\Delta'(\varphi)},
\end{equation}
where the measure is deformed by
\begin{eqnarray}
k^2+{k'}^2&=&1,\\
\Delta(\varphi)&:=&\sqrt{1-k^2\sin^2\varphi},\\
\Delta'(\varphi)&:=&\sqrt{1-{k'}^2\sin^2\varphi}.
\end{eqnarray}

Using these non-circular half angles, and the corresponding ``square
angle'' 
$I_{k'}:=\int_0^{\frac{\pi}2}\frac{d\varphi}{\Delta'(\varphi)}$,  one 
can construct a massive
Dirac spinor wherever the following ``flatness'' condition is
fulfilled:
\begin{eqnarray}
\sum_{a\in\partial F}(I_{k'}-u(a))&=I_{k'} \mod 4I_{k'}&\text{ for 
each 
face }F\in\Lambda_2,\\ 
\sum_{a\ni v}(u(a))&=I_{k'} \mod 4I_{k'}&\text{ for each 
vertex }v\in\Lambda_0.
\end{eqnarray}

Daniel Bennequin
noticed that the fusion algebra of the Ising model could  be
understood at the finite level: Consider a trinion made of cylinders
of a square lattice, of width $m$ and $n$, glued into a cylinder of
width $m+n$. It has been known since Kaufman~\cite{K} that, in the 
transfer
matrix description of the Ising model, the configuration space of the
Ising model on each of the three boundaries is a representation of
spin groups $\spin(m), \spin(n)$ and $\spin(m+n)$ respectively. If $m$
is odd, there exists a unique irreducible representation $\Delta$ of
$\spin(m)$ but when $m$ is even, there are two irreducible
representations, $\Delta^+$ and $\Delta^-$. A pair of pants gives us a
map $\spin(m)\times\spin(n)\to\spin(m+n)$, in the case of a pair of
pants of height zero, it's the inclusion given by the usual
product. The representations of $\spin(m+n)$ induce representations of
the product group that can be split into irreducible
representations. If the three numbers are even,
\begin{eqnarray}
\Delta^+&\to&\Delta^+\otimes\Delta^++\Delta^-\otimes\Delta^-,\\
\Delta^-&\to&\Delta^+\otimes\Delta^-+\Delta^-\otimes\Delta^+,\\
\end{eqnarray}
while if only one of them is even,
\begin{eqnarray}
\Delta&\to&\Delta^+\otimes\Delta+\Delta^-\otimes\Delta,\\
\end{eqnarray}
and if $m$ and $n$ are both odd,
\begin{eqnarray}
\Delta^+&\to&\Delta\otimes\Delta,\\
\Delta^-&\to&\Delta\otimes\Delta.
\end{eqnarray}
Let us compile these data in an array and relabel $\Delta$ by 
$\sigma$, 
$\Delta^+$ by $1$ and $\Delta^-$ by $\epsilon$:
\begin{equation}
\begin{array}{c|ccc}
&1&\epsilon&\sigma\\
\hline
1&1&\epsilon&\sigma\\
\epsilon&\epsilon&1&\sigma\\
\sigma&\sigma&\sigma&1+\epsilon
\end{array}.
\end{equation}
This is read as follows, the $1+\epsilon$ in the slot
$\sigma\otimes\sigma$ for example, means that the representation $1$ 
and the
representation $\epsilon$ of $\spin(m+n)$ both induce a factor
$\sigma\otimes\sigma$ in the representation in the product group
$\spin(m)\times\spin(n)$.

We get exactly the fusion rules of the Ising model. The only
difference compared with the continuous case is that the algebra is
not closed at a finite level. The columns, rows and entries are not
representations of the same group, rather we have a product of
representations of $\spin(n)$ and $\spin(m)$ as a factor of a
representation of $\spin(n+m)$.

These results provide evidence that a discrete conformal field theory
might be looked for: the discrete Dirac spinor at criticality is the
discrete version of the conformal block associated with the field
$\Psi$ and some sort of fusion algebra can be identified at the finite
level. The program we contemplate is, first to investigate other
statistical models and see if there are such patterns. If that is the
case, we must then mimic in the discrete setup the vertex operator
algebra of the continuous conformal theory. This can be attempted by
defining a discrete operator algebra, in a similar fashion to Kadanoff
and Ceva~\cite{KC}, and splitting this algebra according to its
discrete holomorphic and anti-holomorphic parts. The hope is that some
aspects of the powerful results and techniques defined by Belavin,
Polyakov and Zamolodchikov~\cite{BPZ} will still hold. A very
interesting issue would be, as we have done for the Ising model, to
realize the fusion rules of a theory in the discrete setup, yielding
its Verlinde algebra.

\acknowledgements
 The author
did the main part of this work as a Ph.D student in
 Strasbourg, France, under the supervision of Daniel~Bennequin
 \cite{M}; his remarks were essential throughout the paper. The rest
 was done during two postdoctoral stays, in Djursholm, Sweden, funded
 by a Mittag-Leffler Institute grant, and at the University of Tel
 Aviv, Israel, thanks to an Algebraic Lie Representations TMR network
 grant, no. ERB FMRX-CT97-0100. We thank M.~Slupinski, M.~Katz and
 L.~Polterovich for discussions on spinors, M.~Sharir for advice on
 Vorono\"\i~ diagrams, R.~Benedetti for references on flat Riemannian
 metrics, and B.~McCoy for pointing out references in statistical
 mechanics and his great help in understanding the relationship
 between the Ising model and discrete holomorphy.
\endacknowledgements
\bibliographystyle{jcelist} \bibliography{these}

\def\No{\kern-.25em\lower.2ex\hbox{\char'27}}%
\begin{thebibliography}{McCW81}

\bibitem[AdABM]{AdABM}
{ Angl{\`e}s~d'Auriac, J.-Ch., Boukraa, S. and Maillard, J.M.}: Let's
  baxterise. {\em J. Stat. Phys.}, 2000.
\newblock hepth/0003212, To appear

\bibitem[BGB]{BGB}
{ Baake,~M, Grimm,~U and  Baxter,~R.~J.}: A critical {I}sing model on the
  labyrinth. {\em Internat. J. Modern Phys. B}\/ {\bf 8},25-26 (1994),
  3579--3600.
\newblock Perspectives on solvable models

\bibitem[Bax]{Bax}
{ Baxter,~R.~J.}: {\em Exactly solved models in statistical mechanics}. London:
  Academic Press Inc. [Harcourt Brace Jovanovich Publishers],  1989.
\newblock Reprint of the 1982 original

\bibitem[BPZ]{BPZ}
{ Belavin,~A.~A., Polyakov,~A.~M. and Zamolodchikov,~A.~B.}: Infinite conformal
  symmetry in two-dimensional quantum field theory. {\em Nucl. Phys. B}\/
  {\bf 241}, 2, 333--380 (1984)

\bibitem[BS96]{BS96}
{  Benjamini,~I. and Schramm,~O.}: Random walks and harmonic functions on
  infinite planar graphs using square tilings. {\em Ann. Probab.}\/ {\bf 24},3
 , 1219--1238 (1996)

\bibitem[Bus]{Bus}
{  Buser,~P.}: {\em Geometry and spectra of compact {R}iemann surfaces},
   Boston, MA: Birkh\"auser Boston Inc., 1992

\bibitem[CdV90]{CdV90}
{  Colin~de Verdi{\`e}re},~Y.: Un principe variationnel pour les empilements
  de cercles. {\em Invent. Math.}\/ {\bf 104}, 655--669 (1991)

\bibitem[CdV96]{CdV96}
{ Colin~de Verdi{\`e}re,~Y.  Gitler,~I. and Vertigan,~D.}: R\'eseaux
\'electriques planaires. {I}{I}. {\em Comment. Math. Helv.}\/ {\bf
71}, 1, 144--167 (1996)

\bibitem[Cox1]{Cox1}
{ Coxeter,~H.~S.~M.}: {\em Introduction to geometry}. {\em Wiley Classics
  Library},  New York: John Wiley \& Sons Inc., 1989.
\newblock Reprint of the 1969 edition

\bibitem[DD]{DD}
{  Dotsenko,~V.~S. and Dotsenko,~V.~S.}: Critical behaviour of the
  phase transition in the $2$d Ising model with impurities.. {\em Adv. in
  Phys.}\/ {\bf 32}, 2, 129--172 (1983)

\bibitem[Duf]{Duf}
{ Duffin,~R.~J.}: Basic properties of discrete analytic functions. {\em Duke
  Math. J.}\/ {\bf 23}, 335--363 (1956)

\bibitem[GS87]{GS87}
{ Gr{\"u}nbaum,~B. and Shephard,~G.~C.}: {\em Tilings and patterns},
New York: W. H.  Freeman and Company, 1987

\bibitem[Hug]{Hug}
{  Hughes,~B.~D.}: {\em Random walks and random environments. {V}ol. 1},
  New York: The Clarendon Press Oxford University Press, 1995.
\newblock Random walks

\bibitem[ID]{ID}
{ Itzykson,~C. and Drouffe,~J.-M.}: {\em Statistical field theory. 2
{V}ol.}. {\em Cambridge Monographs on Mathematical Physics},
Cambridge: Cambridge University Press, 1989

\bibitem[KC]{KC}
{ Kadanoff,~L.P. and Ceva,~H.}: Determination of an operator algebra
for the two-dimensional {I}sing model. {\em Phys. Rev. B (3)}\/ {\bf
3}, 3918--3939 (1971)

\bibitem[K]{K}
{ Kaufman,~B}: Crystal statistics ii. Partition function evaluated by spinor
  analysis. {\em Phys. Rev.}\/ {\bf 76},1232 (1949)

\bibitem[KO]{KO}
{ Kaufman,~B and Onsager,~L}: Crystal statistics iii. Short range order in a
  binary Ising lattice. {\em Phys. Rev.}\/ {\bf 76},1244 (1949)

\bibitem[Ken]{Ken}
{ Kenyon,~R.}: Tilings and discrete {D}irichlet problems. {\em Israel J.
  Math.}\/ {\bf 105}, 61--84 (1998)

\bibitem[KW]{KW}
{ Kramers,~H.~A. and Wannier,~G.~H.}: Statistics of the two-dimensional
  ferromagnet. {I}. {\em Phys. Rev. (2)}\/ {\bf 60}, 252--262 (1941)

\bibitem[LF]{LF}
{ Lelong-Ferrand,~J.}: {\em Repr\'esentation conforme et
transformations \`a int\'egrale de {D}irichlet born\'ee}. Paris:
Gauthier-Villars, 1955

\bibitem[McCW]{McCW}
{ McCoy,~B.M. and Wu,~T.T}: {\em The two-dimensional Ising model}.
Cambridge, Massachusetts: Harvard University Press, 1973

\bibitem[McCW81]{McCW81}
{ McCoy,~B.M. and Wu,~T.T.}: Non-linear partial difference equations for the
  two-spin correlation function of the two-dimensional Ising model. {\em
  Nucl. Phys. B}\/ {\bf 180}, 89--115 (1981)

\bibitem[M]{M}
{  Mercat,~C.}: {\em Holomorphie discr\`ete et mod\`ele d'Ising}, PhD
  thesis, Universit\'e Louis Pasteur, Strasbourg, France, 1998,
\newblock under the direction of Daniel Bennequin, Pr\'epublication de l'IRMA,
  available at http://www-irma.u-strasbg.fr/irma/publications/1998/98014.shtml

\bibitem[Mil]{Mil}
{ Milnor,~J}: {Spin structures on manifolds}. {\em Enseign. Math., II.
  Ser.}\/ {\bf 9}, 198--203 (1963)

\bibitem[Ons]{Ons}
{  Onsager,~L.}: Crystal statistics. {I}. {A} two-dimensional model with an
  order-disorder transition. {\em Phys. Rev. (2)}\/ {\bf 65}, 117--149 (1944)

\bibitem[PS85]{PS85}
{ Preparata,~F.~P. and  Shamos,~M.~I.}: {\em Computational geometry},
  {\em Texts and Monographs in Computer Science},  New York: Springer-Verlag,
  1985.
\newblock An introduction

\bibitem[SMJ]{SMJ}
{ Sato,~M, Miwa,~T and Jimbo,~M}: Studies on holonomic quantum fields i-iv.
  {\em Proc. Japan Acad.}\/ {\bf 53 (A)},1-6 (1977)

\bibitem[Sie]{Sie}
{ Siegel,~C.~L.}: {\em Topics in complex function theory. {V}ol.\ {I}{I}},
  {\em Wiley Classics Library},  New York: John Wiley \& Sons Inc., 1988
\newblock Automorphic functions and abelian integrals, translated from the
  German by A. Shenitzer and M. Tretkoff, with a preface by Wilhelm Magnus,
  Reprint of the 1971 edition. A Wiley-Interscience Publication

\bibitem[Syo]{Syo}
{  Sy\^ozi, Itiro}: Statistics of kagom\'e lattice. {\em Progress Theoret.
  Phys.}\/ {\bf 6}, 306--308 (1951)

\bibitem[Tro]{Tro}
{ Troyanov, M.}: Les surfaces euclidiennes \`a singularit\'es coniques.
  {\em Enseign. Math. (2)}\/ {\bf 32},1-2, 79--94 (1986)

\bibitem[Uti]{Uti}
{ Utiyama,~T}: Statistics of two-dimensional {I}sing lattices of chequered
  types. {\em Progress Theoret. Phys.}\/ {\bf 6}, 907--909 (1951).
\newblock Letter to the Editor

\bibitem[Veb]{Veb}
{ Veblen, O.}: {\em {Analysis situs. 2. edit.}}, New York: American
Mathematical Society X, 1931

\bibitem[VoSh]{VoSh}
{ Voevodski{\u\i},~V.~A. and Shabat,~G.~B.}: Equilateral triangulations of
  {R}iemann surfaces, and curves over algebraic number fields. {\em Dokl. Akad.
  Nauk SSSR}\/ {\bf 304}, 2, 265--268  (1989)

\bibitem[Wan50]{Wan50}
{ Wannier,~G.~H.}: Antiferromagnetism. {T}he triangular {I}sing net. {\em
  Physical Rev. (2)}\/ {\bf 79}, 357--364 (1950)

\bibitem[Whit]{Whit}
{ Whitney,~H}: Product on complexes.. {\em Annals of Math. (2)}\/ {\bf 39}
 , 397--432 (1938)

\bibitem[Yam]{Yam}
{  Yamamoto, Tunenobu}: On the crystal statistics of two-dimensional {I}sing
  ferromagnets. {\em Progress Theoret. Phys.}\/ {\bf 6}, 533--542 (1951)

\end{thebibliography}

\end{document}